\def\year{2016}\relax
\documentclass[letterpaper]{article}
\usepackage{aaai16}
\usepackage{times}
\usepackage{helvet}
\usepackage{courier}
\usepackage{url}
\usepackage{mdwmath}
\usepackage{mathrsfs}
\usepackage{multirow}
\usepackage{mdwtab}
\usepackage{amssymb}
\usepackage{amsthm}
\usepackage{amsmath}
\usepackage{array}
\usepackage{amsfonts}
\usepackage{algorithm}
\usepackage{algorithmic}
\usepackage[tight,footnotesize]{subfigure}
\usepackage{eqparbox}
\usepackage{enumerate}
\usepackage{framed}
\usepackage{fancyhdr}
\usepackage{graphicx}
\usepackage{lipsum}
\usepackage{enumitem}

\newcommand{\citet}[1]{\citeauthor{#1}~\shortcite{#1}}
\newcommand{\citep}{\cite}

\newtheorem{definition}{Definition}
\newtheorem{theorem}{Theorem}
\newtheorem{lemma}{Lemma}

\newtheorem{property}{Property}

\begin{document}

\title{Scalable Algorithms for Tractable Schatten Quasi-Norm Minimization}

\author{Fanhua Shang, Yuanyuan Liu, James Cheng\\
Department of Computer Science and Engineering, The Chinese University of Hong Kong\\
\{fhshang, yyliu, jcheng\}@cse.cuhk.edu.hk
}
\maketitle

\begin{abstract}
The Schatten-$p$ quasi-norm $(0\!\!<\!\!p\!\!<\!\!1)$ is usually used to replace the standard nuclear norm in order to approximate the rank function more accurately. However, existing Schatten-$p$ quasi-norm minimization algorithms involve singular value decomposition (SVD) or eigenvalue decomposition (EVD) in each iteration, and thus may become very slow and impractical for large-scale problems. In this paper, we first define two tractable Schatten quasi-norms, i.e., the Frobenius/nuclear hybrid and bi-nuclear quasi-norms, and then prove that they are in essence the Schatten-$2/3$ and $1/2$ quasi-norms, respectively, which lead to the design of very efficient algorithms that only need to update two much smaller factor matrices. We also design two efficient proximal alternating linearized minimization algorithms for solving representative matrix completion problems. Finally, we provide the global convergence and performance guarantees for our algorithms, which have better convergence properties than existing algorithms. Experimental results on synthetic and real-world data show that our algorithms are more accurate than the state-of-the-art methods, and are orders of magnitude faster.
\end{abstract}

\section{Introduction}
\noindent In recent years, the matrix rank minimization problem arises in a wide range of applications such as matrix completion, robust principal component analysis, low-rank representation, multivariate regression and multi-task learning. To solve such problems, \citet{fazel:rmh,candes:crmc,recht:nnm} have suggested to relax the rank function by its convex envelope, i.e., the nuclear norm. In fact, the nuclear norm is equivalent to the $\ell_{1}$-norm on singular values of a matrix, and thus it promotes a low-rank solution. However, it has been shown in~\citep{fan:ve} that the $\ell_{1}$-norm regularization over-penalizes large entries of vectors, and results in a biased solution. By realizing the intimate relationship between them, the nuclear norm penalty also over-penalizes large singular values, that is, it may make the solution deviate from the original solution as the $\ell_{1}$-norm does~\citep{nie:lrmr,lu:gsvt}. Compared with the nuclear norm, the Schatten-${p}$ quasi-norm for $0\!\!<\!\!p\!\!<\!\!1$ makes a closer approximation to the rank function. Consequently, the Schatten-${p}$ quasi-norm minimization has attracted a great deal of attention in images recovery~\citep{lu:irsvm,lu:lrm}, collaborative filtering~\citep{nie:rmc,lu:gsvt,mohan:mrm} and MRI analysis~\citep{majumdar:mri}. In addition, many non-convex surrogate functions of the $\ell_{0}$-norm listed in~\citep{lu:lrm,lu:gsvt} have been extended to approximate the rank function, such as SCAD~\citep{fan:ve} and MCP~\citep{zhang:mcp}.

All non-convex surrogate functions mentioned above for low-rank minimization lead to some non-convex, non-smooth, even non-Lipschitz optimization problems. Therefore, it is crucial to develop fast and scalable algorithms which are specialized to solve some alternative formulations. So far, \citet{lai:irls} proposed an iterative reweighted lease squares (IRucLq) algorithm to approximate the Schatten-${p}$ quasi-norm minimization problem, and proved that the limit point of any convergent subsequence generated by their algorithm is a critical point. Moreover, \citet{lu:lrm} proposed an iteratively reweighted nuclear norm (IRNN) algorithm to solve many non-convex surrogate minimization problems. For matrix completion problems, the Schatten-${p}$ quasi-norm has been shown to be empirically superior to the nuclear norm~\citep{marjanovic:mc}. In addition, \citet{zhang:ncmr} theoretically proved that the Schatten-${p}$ quasi-norm minimization with small $p$ requires significantly fewer measurements than the convex nuclear norm minimization. However, all existing algorithms have to be solved iteratively and involve SVD or EVD in each iteration, which incurs high computational cost and is too expensive for solving large-scale problems~\citep{cai:fsvt,liu:nnr}.

In contrast, as an alternative non-convex formulation of the nuclear norm, the bilinear spectral regularization as in~\citep{srebro:mmmf,recht:nnm} has been successfully applied in many large-scale applications, e.g., collaborative filtering~\citep{mitra:lsmf}. As the Schatten-${p}$ quasi-norm is equivalent to the $\ell_{p}$ quasi-norm on singular values of a matrix, it is natural to ask the following question: \emph{can we design equivalent matrix factorization forms for the cases of the Schatten quasi-norm, e.g., $p=2/3$ or $1/2$}?

In order to answer the above question, in this paper we first define two tractable Schatten quasi-norms, i.e., the Frobenius/nuclear hybrid and bi-nuclear quasi-norms. We then prove that they are in essence the Schatten-${2/3}$ and ${1/2}$ quasi-norms, respectively, for solving whose minimization we only need to perform SVDs on two much smaller factor matrices as contrary to the larger ones used in existing algorithms, e.g., IRNN. Therefore, our method is particularly useful for many ``big data'' applications that need to deal with large, high dimensional data with missing values. To the best of our knowledge, this is the first paper to scale Schatten quasi-norm solvers to the Netflix dataset. Moreover, we provide the global convergence and recovery performance guarantees for our algorithms. In other words, this is the best guaranteed convergence for algorithms that solve such challenging problems.

\section{Notations and Background}
The Schatten-${p}$ norm ($0<p<\infty$) of a matrix $X\in \mathbb{R}^{m\times n}$ ($m\geq n$) is defined as
\vspace{-2mm}
\begin{displaymath}
\|X\|_{S_{p}}\triangleq\left(\sum^{n}_{i=1}\sigma^{p}_{i}(X)\right)^{1/p},
\end{displaymath}
where $\sigma_{i}(X)$ denotes the $i$-th singular value of $X$. When $p\!=\!1$, the Schatten-${1}$ norm is the well-known nuclear norm, $\|X\|_{*}$. In addition, as the non-convex surrogate for the rank function, the Schatten-${p}$ quasi-norm with $0\!<\!p\!<\!1$ is a better approximation than the nuclear norm~\citep{zhang:ncmr} (analogous to the superiority of the $\ell_{p}$ quasi-norm to the $\ell_{1}$-norm~\citep{daubechies:rlsm}).

To recover a low-rank matrix from some linear observations $b\!\in\!\mathbb{R}^{s}$, we consider the following general Schatten-${p}$ quasi-norm minimization problem,
\vspace{-1mm}
\begin{equation}\label{Nb1}
\min_{X} \lambda\|X\|^{p}_{S_{p}}+f\!\left(\mathcal{A}(X)-b\right),
\end{equation}
where $\mathcal{A}\!:\!\mathbb{R}^{m\times n}\!\rightarrow\! \mathbb{R}^{s}$ denotes the linear measurement operator, $\lambda\!>\!0$ is a regularization parameter, and the loss function $f(\cdot):\mathbb{R}^{s}\!\rightarrow\! \mathbb{R}$ generally denotes certain measurement for characterizing $\mathcal{A}(X)-b$. The above formulation can address a wide range of problems, such as matrix completion~\citep{marjanovic:mc,rohde:lrm} ($\mathcal{A}$ is the sampling operator and $f(\cdot)\!=\!\|\!\cdot\!\|^{2}_{2}$), robust principal component analysis~\citep{candes:rpca,wang:rmr,shang:rpca} ($\mathcal{A}$ is the identity operator and $f(\cdot)\!=\!\|\!\cdot\!\|_{1}$), and multivariate regression~\citep{hsieh:nnm} ($\mathcal{A}(X)\!=\!AX$ with $A$ being a given matrix, and $f(\cdot)\!=\!\|\!\cdot\!\|^{2}_{F}$). Furthermore, $f(\cdot)$ may be also chosen as the Hinge loss in~\citep{srebro:mmmf} or the $\ell_{p}$ quasi-norm in~\citep{nie:rmc}.

Analogous to the $\ell_{p}$ quasi-norm, the Schatten-${p}$ quasi-norm is also non-convex for $p\!<\!1$, and its minimization is generally NP-hard~\citep{lai:irls}.
Therefore, it is crucial to develop efficient algorithms to solve some alternative formulations of Schatten-${p}$ quasi-norm minimization \eqref{Nb1}. So far, only few algorithms~\citep{lai:irls,mohan:mrm,nie:rmc,lu:lrm} have been developed to solve such problems. Furthermore, since all existing Schatten-${p}$ quasi-norm minimization algorithms involve SVD or EVD in each iteration, they suffer from a high computational cost of $O(n^{2}m)$, which severely limits their applicability to large-scale problems. Although there have been many efforts towards fast SVD or EVD computation such as partial SVD~\citep{larsen:svd}, the performance of those methods is still unsatisfactory for real-life applications~\citep{cai:fsvt}.

\section{Tractable Schatten Quasi-Norms}
As in~\citep{srebro:mmmf}, the nuclear norm has the following alternative non-convex formulations.
\begin{lemma}
Given a matrix $X\in \mathbb{R}^{m\times n}$ with $\textrm{rank}(X)=r\leq d$, the following holds:
\vspace{-2mm}
\begin{equation*}
\begin{split}
\|X\|_{*}&=\min_{U\in\mathbb{R}^{m\!\times\! d},V\in\mathbb{R}^{n\!\times\! d}:X=UV^{T}}\|U\|_{F}\|V\|_{F}\\
&=\min_{U,V:X=UV^{T}}\frac{\|U\|^{2}_{F}+\|V\|^{2}_{F}}{2}.
\end{split}
\end{equation*}
\end{lemma}

\subsection{Frobenius/Nuclear Hybrid Quasi-Norm}
Motivated by the equivalence relation between the nuclear norm and the bilinear spectral regularization (please refer to~\citep{srebro:mmmf,recht:nnm}), we define a Frobenius/nuclear hybrid (F/N) norm as follows

\begin{definition}
For any matrix $X\in \mathbb{R}^{m\times n}$ with $\textrm{rank}(X)=r\leq d$, we can factorize it into two much smaller matrices $U\in \mathbb{R}^{m\times d}$ and $V\in \mathbb{R}^{n\times d}$ such that $X=UV^{T}$. Then the Frobenius/nuclear hybrid norm of $X$ is defined as
\begin{equation*}
\|X\|_{\textup{F/N}}:=\min_{X=UV^{T}}\|U\|_{*}\|V\|_{F}.
\end{equation*}
\end{definition}
\vspace{-1mm}

In fact, the Frobenius/nuclear hybrid norm is not a real norm, because it is non-convex and does not satisfy the triangle inequality of a norm. Similar to the well-known Schatten-${p}$ quasi-norm ($0\!<\!p\!<\!1$), the Frobenius/nuclear hybrid norm is also a quasi-norm, and their relationship is stated in the following theorem.

\begin{theorem}
The Frobenius/nuclear hybrid norm $\|\!\cdot\!\|_{\textup{F/N}}$ is a quasi-norm. Surprisingly, it is also the Schatten-${2/3}$ quasi-norm, i.e.,
\begin{equation*}
\|X\|_{\textup{F/N}}=\|X\|_{S_{2/3}},
 \end{equation*}
where $\|X\|_{S_{2/3}}$ denotes the Schatten-${2/3}$ quasi-norm of $X$.
\end{theorem}

\begin{property}
For any matrix $X\in \mathbb{R}^{m\times n}$ with $\textrm{rank}(X)=r\leq d$, the following holds:
\begin{equation*}
\begin{split}
\|X\|_{\textup{F/N}}=&\min_{U\in\mathbb{R}^{m\!\times\! d},V\in\mathbb{R}^{n\!\times\! d}:X=UV^{T}}\|U\|_{*}\|V\|_{F}\\
=&\min_{X=UV^{T}}\!\left(\frac{2\|U\|_{*}+\|V\|^{2}_{F}}{3}\right)^{3/2}.
\end{split}
\end{equation*}
\end{property}

The proofs of Property 1 and Theorem 1 can be found in the Supplementary Materials.

\subsection{Bi-Nuclear Quasi-Norm}
Similar to the definition of the above Frobenius/nuclear hybrid norm, our bi-nuclear (BiN) norm is naturally defined as follows.

\begin{definition}
For any matrix $X\in \mathbb{R}^{m\times n}$ with $\textrm{rank}(X)=r\leq d$, we can factorize it into two much smaller matrices $U\in \mathbb{R}^{m\times d}$ and $V\in \mathbb{R}^{n\times d}$ such that $X=UV^{T}$. Then the bi-nuclear norm of $X$ is defined as
\begin{equation*}
\|X\|_{\textup{BiN}}:=\min_{X=UV^{T}}\|U\|_{*}\|V\|_{*}.
\end{equation*}
\end{definition}

Similar to the Frobenius/nuclear hybrid norm, the bi-nuclear norm is also a quasi-norm, as stated in the following theorem.
\begin{theorem}
The bi-nuclear norm $\|\!\cdot\!\|_{\textup{BiN}}$ is a quasi-norm. In addition, it is also the Schatten-${1/2}$ quasi-norm, i.e.,
\begin{equation*}
\|X\|_{\textup{BiN}}=\|X\|_{S_{1/2}}.
 \end{equation*}
\end{theorem}

The proof of Theorem 2 can be found in the Supplementary Materials. Due to the relationship between the bi-nuclear quasi-norm and the Schatten-${1/2}$ quasi-norm, it is easy to verify that the bi-nuclear quasi-norm possesses the following properties.

\begin{property}
For any matrix $X\in \mathbb{R}^{m\times n}$ with $\textrm{rank}(X)=r\leq d$, the following holds:
\vspace{-2mm}
\begin{equation*}
\begin{split}
\|X\|_{\textup{BiN}}&=\!\min_{X=UV^{T}}\!\|U\|_{*}\|V\|_{*}\!=\!\min_{X=UV^{T}}\!\frac{\|U\|^{2}_{*}\!+\!\|V\|^{2}_{*}}{2}\\
&=\min_{X=UV^{T}}\!\left(\frac{\|U\|_{*}\!+\!\|V\|_{*}}{2}\right)^{2}.
\end{split}
\end{equation*}
\end{property}

The following relationship between the nuclear norm and the Frobenius norm is well known: $\|X\|_{F}\!\leq\!\|X\|_{*}\!\leq\! \sqrt{r}\|X\|_{F}$. Similarly, the analogous bounds hold for the Frobenius/nuclear hybrid and bi-nuclear quasi-norms, as stated in the following property.

\begin{property}
For any matrix $X\in \mathbb{R}^{m\times n}$ with $\textrm{rank}(X)= r$, the following inequalities hold:
\begin{displaymath}
\begin{split}
\|X\|_{*}\leq\|X\|_{\textup{F/N}} & \leq\sqrt{r}\|X\|_{*},\\
\|X\|_{*}\leq\|X\|_{\textup{F/N}}\leq & \|X\|_{\textup{BiN}}\leq r\|X\|_{*}.
\end{split}
\end{displaymath}
\end{property}
The proof of Property 3 can be found in the Supplementary Materials. It is easy to see that Property 3 in turn implies that any low Frobenius/nuclear hybrid or bi-nuclear norm approximation is also a low nuclear norm approximation.

\section{Optimization Algorithms}

\subsection{Problem Formulations}
To bound the Schatten-${2/3}$ or -${1/2}$ quasi-norm of $X$ by $\frac{1}{3}(2\|U\|_{*}\!+\!\|V\|^{2}_{F})$ or $\frac{1}{2}(\|U\|_{*}\!+\!\|V\|_{*})$, we mainly consider the following general structured matrix factorization formulation as in~\citep{haeffele:lrmf},
\begin{equation}\label{Fm1}
\min_{U,V}\;\lambda\varphi(U,V)+f(\mathcal{A}(UV^{T})-b),
\end{equation}
where the regularization term $\varphi(U,V)$ denotes $\frac{1}{3}(2\|U\|_{*}\!+\!\|V\|^{2}_{F})$ or $\frac{1}{2}(\|U\|_{*}\!+\!\|V\|_{*})$.

As mentioned above, there are many Schatten-${p}$ quasi-norm minimization problems for various real-world applications. Therefore, we propose two efficient algorithms to solve the following low-rank matrix completion problems:
\begin{equation}\label{Fm2}
\min_{U,V} \frac{\lambda(2\|U\|_{*}\!+\!\|V\|^{2}_{F})}{3}+\frac{1}{2}\|\mathcal{P}_{\Omega}(UV^{T})\!-\!\mathcal{P}_{\Omega}(D)\|^{2}_{F},
\end{equation}
\begin{equation}\label{Fm3}
\min_{U,V} \frac{\lambda(\|U\|_{*}+\|V\|_{*})}{2}+\frac{1}{2}\|\mathcal{P}_{\Omega}(UV^{T})-\mathcal{P}_{\Omega}(D)\|^{2}_{F},
\end{equation}
where $\mathcal{P}_{\Omega}$ denotes the linear projection operator, i.e., $\mathcal{P}_{\Omega}(D)_{ij}\!\!=\!\!D_{ij}$ if $(i,j)\!\!\in\!\!\Omega$, and $\mathcal{P}_{\Omega}(D)_{ij}\!\!=\!\!0$ otherwise. Due to the operator $\mathcal{P}_{\Omega}$ in \eqref{Fm2} and \eqref{Fm3}, we usually need to introduce some auxiliary variables for solving them. To avoid introducing auxiliary variables, motivated by the proximal alternating linearized minimization (PALM) method proposed in~\citep{bolte:palm}, we propose two fast PALM algorithms to efficiently solve \eqref{Fm2} and \eqref{Fm3}. The space limitation refrains us from fully describing each algorithm, but we try to give enough details of a representative algorithm for solving \eqref{Fm2} and discussing their differences.

\subsection{Updating $U_{k+1}$ and $V_{k+1}$ with Linearization Techniques}
Let $g_{k}(U)\!:=\!\|\mathcal{P}_{\Omega}(UV^{T}_{k})\!-\!\mathcal{P}_{\Omega}(D)\|^{2}_{F}/2$, and then its gradient is Lipschitz continuous with constant $l^{g}_{k+1}$, meaning that $\|\nabla\! g_{k}(U_{1})\!-\!\nabla\! g_{k}(U_{2})\|_{F}\!\leq\! l^{g}_{k+1}\|U_{1}\!-\!U_{2}\|_{F}$ for any $U_{1}, U_{2}\!\in\! \mathbb{R}^{m\times d}$. By linearizing $g_{k}(U)$ at $U_{k}$ and adding a proximal term, then we have the following approximation:
\vspace{-2mm}
\begin{equation}\label{Al1}
\widehat{g_{k}}(U,U_{k}\!)\!=\!g_{k}(U_{k}\!)\!+\!\langle\nabla\! g_{k}(U_{k}\!),U\!-\!U_{k}\rangle\!+\!\frac{l^{g}_{k\!+\!1}}{2}\|U\!-\!U_{k}\!\|^{2}_{F}.
\end{equation}
Thus, we have
\vspace{-2mm}
\begin{equation}\label{Al2}
\begin{split}
&U_{k+1}\!=\!\mathop{\arg\min}_{U}\frac{2\lambda}{3}\|U\|_{*}\!+\!\widehat{g_{k}}(U,U_{k})\\
=&\!\mathop{\arg\min}_{U}\!\frac{2\lambda}{3}\|U\|_{*}\!+\!\frac{l^{g}_{k\!+\!1}}{2}\|U\!-\!U_{k}\!+\!\frac{\nabla\! g_{k}(U_{k})}{l^{g}_{k\!+\!1}}\|^{2}_{F}.
\end{split}
\end{equation}

Similarly, we have
\vspace{-2mm}
\begin{equation}\label{Al3}
V_{k\!+\!1}\!\!=\!\mathop{\arg\min}_{V}\!\!\frac{\lambda}{3}\|V\|^{2}_{F}\!+\!\frac{l^{h}_{k\!+\!1}}{2}\|V\!-\!V_{k}\!+\!\nabla\! h_{k}(V_{k})/l^{h}_{k\!+\!1}\|^{2}_{F},
\end{equation}
where $h_{k}(V):=\|\mathcal{P}_{\Omega}(U_{k+1}V^{T})\!-\!\mathcal{P}_{\Omega}(D)\|^{2}_{F}/2$ with the Lipschitz constant $l^{h}_{k+1}$. The problems \eqref{Al2} and \eqref{Al3} are known to have closed-form solutions, which of the former is given by the so-called matrix shrinkage operator~\citep{cai:svt}. In contrast, for solving \eqref{Fm3}, $U_{k+1}$ is computed in the same way as \eqref{Al2}, and $V_{k+1}$ is given by
\vspace{-2mm}
\begin{equation}\label{Al4}
V_{k+1}\!\!=\!\mathop{\arg\min}_{V}\!\!\frac{\lambda}{2}\|V\|_{*}\!+\!\frac{l^{h}_{k\!+\!1}}{2}\|V\!-\!V_{k}\!+\!\nabla\! h_{k}(V_{k})/l^{h}_{k\!+\!1}\|^{2}_{F}.
\end{equation}

\subsection{Updating Lipschitz Constants}
Next we compute the Lipschitz constants $l^{g}_{k+1}$ and $l^{h}_{k+1}$ at the $(k\!+\!1)$-iteration.
\vspace{-1mm}
\begin{equation*}
\begin{split}
&\|\nabla g_{k}(U_{1})\!-\!\nabla g_{k}(U_{2})\|_{F}\!=\!\|[\mathcal{P}_{\Omega}(U_{1}V^{T}_{k}-U_{2}V^{T}_{k})]V_{k}\|_{F}\\
\leq & \|V_{k}\|^{2}_{2}\|U_{1}-U_{2}\|_{F},\\
&\!\!\|\nabla h_{k}(V_{1})\!-\!\nabla h_{k}(V_{2})\|_{F}\!=\!\|U^{T}_{k+1}\![\mathcal{P}_{\Omega}(U_{k+1}\!(V^{T}_{1}\!-\!V^{T}_{2}))]\|_{F}\\
\leq& \|U_{k+1}\|^{2}_{2}\|V_{1}\!-\!V_{2}\|_{F}.
\end{split}
\end{equation*}
Hence, both Lipschitz constants are updated by
\begin{equation}\label{Al5}
l^{g}_{k+1}=\|V_{k}\|^{2}_{2}\;\,\textrm{and}\;\,l^{h}_{k+1}=\|U_{k+1}\|^{2}_{2}.
\end{equation}

\subsection{PALM Algorithms}
Based on the above development, our algorithm for solving \eqref{Fm2} is given in Algorithm 1. Similarly, we also design an efficient PALM algorithm for solving \eqref{Fm3}. The running time of Algorithm 1 is dominated by performing matrix multiplications. The total time complexity of Algorithm 1, as well as the algorithm for solving \eqref{Fm3}, is $O(nmd)$, where $d\ll m,n$.

\begin{algorithm}[t]
\caption{Solving \eqref{Fm2} via PALM}
\label{alg:Framwork2}
\renewcommand{\algorithmicrequire}{\textbf{Input:}}
\renewcommand{\algorithmicensure}{\textbf{Initialize:}}
\renewcommand{\algorithmicoutput}{\textbf{Output:}}
\begin{algorithmic}[1]
\REQUIRE $\mathcal{P}_{\Omega}(D)$, the given rank $d$ and $\lambda$.
\ENSURE $U_{0}$, $V_{0}$, $\varepsilon$ and $k=0$.\\
\WHILE {not converged}
\STATE {Update $l^{g}_{k+1}$ and $U_{k+1}$ by \eqref{Al5} and \eqref{Al2}, respectively.}
\STATE {Update $l^{h}_{k+1}$ and $V_{k+1}$ by \eqref{Al5} and \eqref{Al3}, respectively.}
\STATE {Check the convergence condition,\\
$\quad \max\{\|U_{k+1}\!-\!U_{k}\|_{F}, \|V_{k+1}\!-\!V_{k}\|_{F}\}<\varepsilon$.}
\ENDWHILE
\OUTPUT $U_{k+1}$, $V_{k+1}$.
\end{algorithmic}
\end{algorithm}

\section{Algorithm Analysis}
We now provide the global convergence and low-rank matrix recovery guarantees for Algorithm 1, and the similar results can be obtained for the algorithm for solving (4).

\subsection{Global Convergence}
Before analyzing the global convergence of Algorithm 1, we first introduce the definition of the critical points of a non-convex function given in~\citep{bolte:palm}.

\begin{definition}
Let a non-convex function $f\!:\!\mathbb{R}^{n}\!\rightarrow\!(-\infty, +\infty]$ be a proper and lower semi-continuous function, and $\textup{dom}f\!=\!\{x\in\mathbb{R}^{n}:f(x)<+\infty\}$.
\begin{itemize}
\item For any $x\in\textup{dom}f$, the Fr\`{e}chet sub-differential of $f$ at $x$ is defined as
\vspace{-2mm}
  \begin{displaymath}
  \widehat{\partial}f(x)\!=\!\{u\!\in\!\mathbb{R}^{n}\!:\lim_{y\neq x}\inf_{y\rightarrow x}\!\frac{f(y)\!-\!f(x)\!-\!\langle u,y\!-\!x\rangle}{\|y\!-\!x\|_{2}}\!\geq\!0\},
  \end{displaymath}
  and $\widehat{\partial}f(x)=\emptyset$ if $x\notin\textup{dom}f$.
\item The limiting sub-differential of $f$ at $x$ is defined as
\vspace{-1mm}
  \begin{displaymath}
  \begin{split}
  \partial f(x)\!=\!\{u\!\in\!\mathbb{R}^{n}\!:\exists x^{k}\rightarrow x,\; f(x^{k})\rightarrow f(x)\\
  \textup{and}\;\,u^{k}\!\in\!\widehat{\partial}f(x^{k})\!\rightarrow\! u\;\textup{as}\;k\!\rightarrow\!\infty\}.
  \end{split}
  \end{displaymath}
\item The points whose sub-differential contains $0$ are called critical points. For instance, the point $x$ is a critical point of $f$ if $0\!\in\!\partial f(x)$.
\end{itemize}
\end{definition}

\begin{theorem} [Global Convergence]
Let $\{(U_{k},V_{k})\}$ be a sequence generated by Algorithm 1, then it is a Cauchy sequence and converges to a critical point of \eqref{Fm2}.
\end{theorem}

The proof of the theorem can be found in the Supplementary Materials. Theorem 3 shows the global convergence of Algorithm 1. We emphasize that, different from the general subsequence convergence property, the global convergence property is given by $(U_{k},V_{k})\!\rightarrow\!(\widehat{U},\widehat{V})$ as the number of iteration $k\!\rightarrow\!+\infty$, where $(\widehat{U},\widehat{V})$ is a critical point of \eqref{Fm2}. As we have stated, existing algorithms for solving the non-convex and non-smooth problem, such as IRucLq and IRNN, have only subsequence convergence~\citep{xu:bcd}. According to~\citep{attouch:cr}, we know that the convergence rate of Algorithm 1 is at least sub-linear, as stated in the following theorem.

\begin{theorem} [Convergence Rate]
The sequence $\{(U_{k},V_{k})\}$ generated by Algorithm 1 converges to a critical point $(\widehat{U},\widehat{V})$ of \eqref{Fm2} at least in the sub-linear convergence rate, that is, there exists $C>0$ and $\theta\in(1/2,1)$ such that
\vspace{-2mm}
\begin{displaymath}
\|[U^{T}_{k},V^{T}_{k}]-[\widehat{U}^{T},\widehat{V}^{T}]\|_{F}\leq Ck^{-\frac{1-\theta}{2\theta-1}}.
\end{displaymath}
\end{theorem}

\subsection{Recovery Guarantee}
In the following, we show that when sufficiently many entries are observed, the critical point generated by our algorithms recovers a low-rank matrix ``close to" the ground-truth one. Without loss of generality, assume that $D\!=\!Z\!+\!E\!\in \!\mathbb{R}^{m\times n}$, where $Z$ is a true matrix, and $E$ denotes a random gaussian noise.

\begin{theorem}
Let $(\widehat{U},\widehat{V})$ be a critical point of the problem \eqref{Fm2} with given rank $d$, and $m\geq n$. Then there exists an absolute constant $C_{1}$, such that with probability at least $1-2\exp(-m)$,
\vspace{-2mm}
\begin{equation*}
\frac{\|Z\!-\!\widehat{U}\widehat{V}^{T}\!\|_{F}}{\sqrt{m n}}\!\!\leq\!\! \frac{\|E\|_{F}}{\sqrt{mn}}\!+\!C_{1}\beta\!\left(\!\frac{md\log(m)}{|\Omega|}\!\right)^{1/4}\!+\!\frac{2\sqrt{d}\lambda}{3C_{2}\sqrt{|\Omega|}},
\end{equation*}
where $\beta=\max_{i,j}|D_{i,j}|$ and $C_{2}=\frac{\|\mathcal{P}_{\Omega}(D-\hat{U}\hat{V}^{T})\hat{V}\|_{F}}{\|\mathcal{P}_{\Omega}(D-\hat{U}\hat{V}^{T})\|_{F}}$.
\end{theorem}

The proof of the theorem and the analysis of lower-boundedness of $C_{2}$ can be found in the Supplementary Materials. When the samples size $|\Omega|\!\gg\! md\log(m)$, the second and third terms diminish, and the recovery error is essentially bounded by the ``average" magnitude of entries of the noise matrix $E$. In other words, only $O(md\log(m))$ observed entries are needed, which is significantly lower than $O(mr\log^{2}(m))$ in standard matrix completion theories~\citep{candes:emc,keshavan:mc,recht:mc}. We will confirm this result by our experiments in the following section.

\begin{figure*}[!t]
\begin{center}
\subfigure[20\% SR and $nf\!=\!0.1$]{\includegraphics[width=0.497\columnwidth]{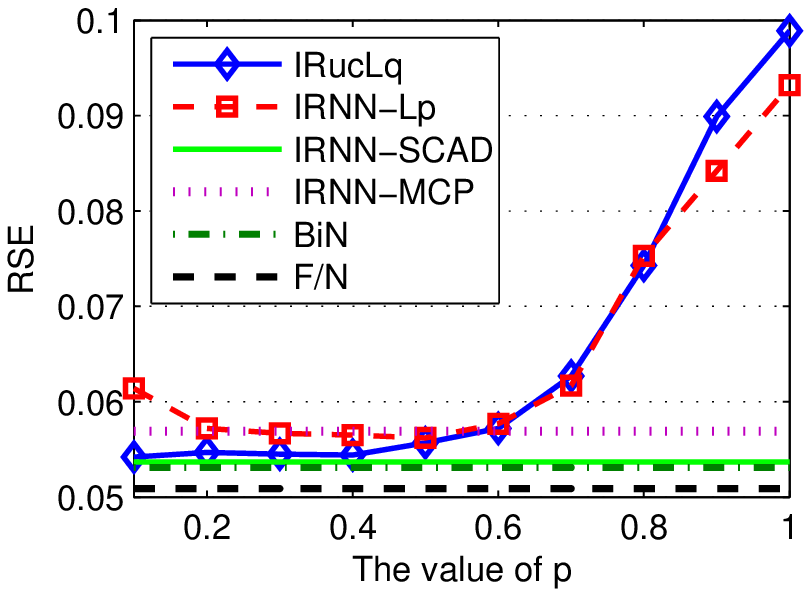}}
\subfigure[20\% SR and $nf\!=\!0.2$]{\includegraphics[width=0.497\columnwidth]{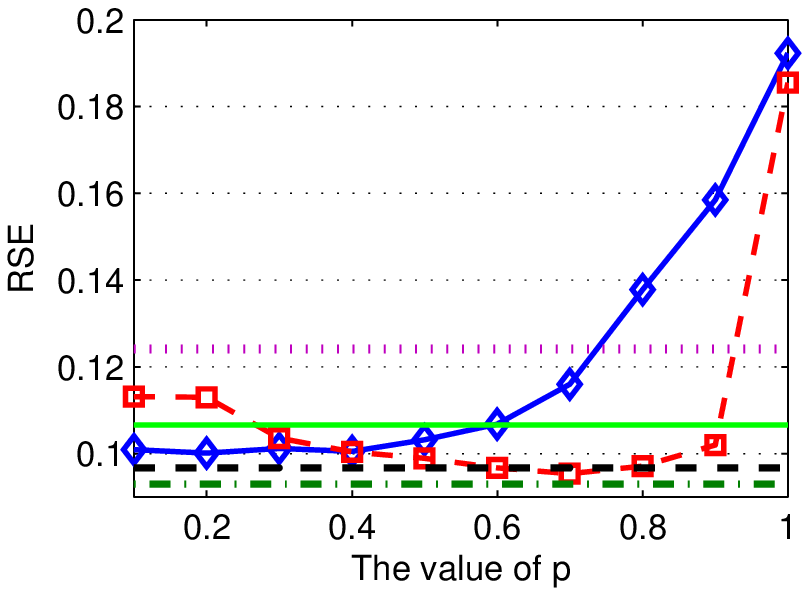}}
\subfigure[30\% SR and $nf\!=\!0.1$]{\includegraphics[width=0.497\columnwidth]{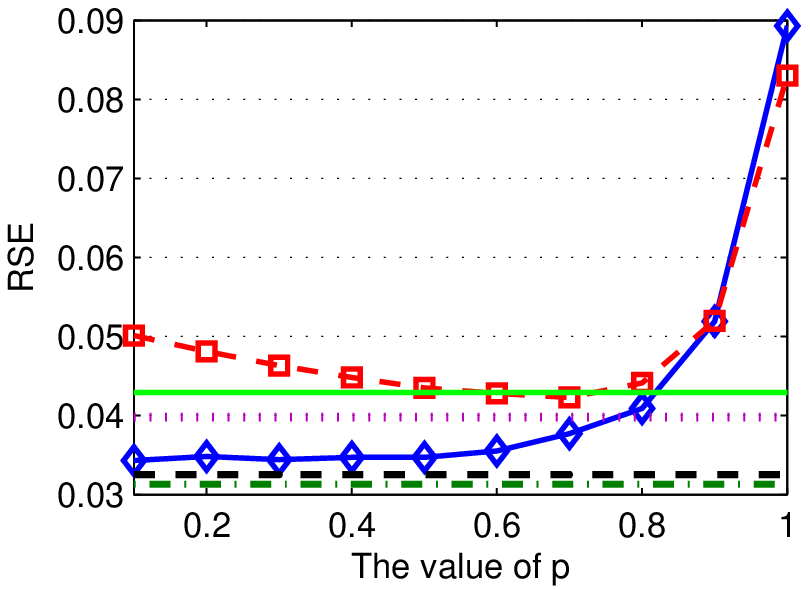}}
\subfigure[30\% SR and $nf\!=\!0.2$]{\includegraphics[width=0.499\columnwidth]{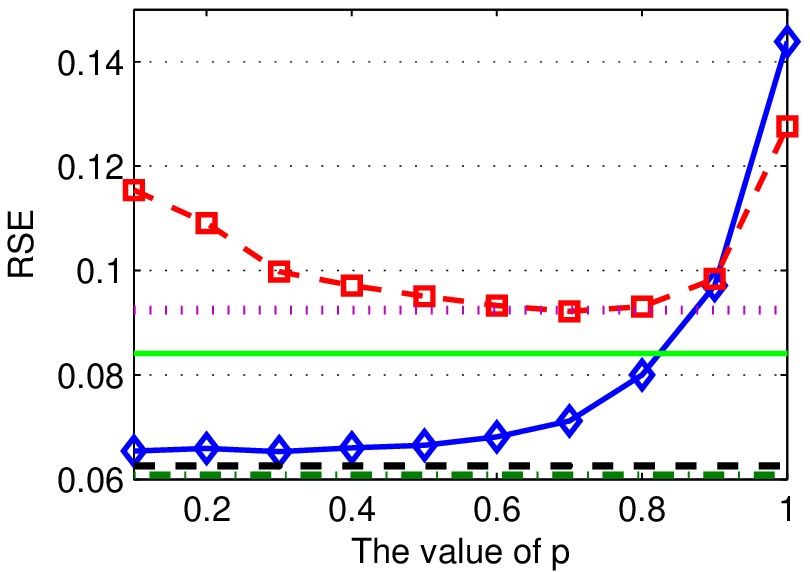}}
\caption{The recovery accuracy of IRucLq, IRNN and our algorithms on noisy random matrices of size $100\times100$.}
\label{fig_sim1}
\end{center}
\end{figure*}

\section{Experimental Results}
We now evaluate both the effectiveness and efficiency of our algorithms for solving matrix completion problems, such as collaborative filtering and image recovery. All experiments were conducted on an Intel Xeon E7-4830V2 2.20GHz CPU with 64G RAM.

\subsubsection{Algorithms for Comparison}
We compared our algorithms, BiN and F/N, with the following state-of-the-art methods: IRucLq{\footnote{\url{http://www.math.ucla.edu/~wotaoyin/}}}~\citep{lai:irls}: In IRucLq, $p$ varies from $0.1$ to $1$ with increment 0.1, and the parameters $\lambda$ and $\alpha$ are set to $10^{-6}$ and $0.9$, respectively. In addition, the rank parameter of the algorithm is updated dynamically as in~\citep{lai:irls}, that is, it only needs to compute the partial EVD. IRNN{\footnote{\url{https://sites.google.com/site/canyilu/}}}~\citep{lu:lrm}: We choose the $\ell_{p}$-norm, SCAD and MCP penalties as the regularization term among eight non-convex penalty functions, where $p$ is chosen from the range of $\{0.1,0.2,\ldots,1\}$. At each iteration, the parameter $\lambda$ is dynamically decreased by $\lambda_{k}\!=\!0.7\lambda_{k-1}$, where $\lambda_{0}\!=\!10\|\mathcal{P}_{\Omega}(D)\|_{\infty}$.

For our algorithms, we set the regularization parameter $\lambda\!=\!5$ or $\lambda\!=\!100$ for noisy synthetic and real-world data, respectively. Note that the rank parameter $d$ is estimated by the strategy in~\citep{wen:nsor}. In addition, we evaluate the performance of matrix recovery by the relative squared error (RSE) and the root mean square error (RMSE), i.e., $\textup{RSE}\!:=\!\|X\!-\!Z\|_{F}/\|Z\|_{F}$ and $\textup{RMSE}\!:=\!\frac{1}{|T|}\!\sqrt{\Sigma_{(i,j)\in T}(X_{ij}\!-\!D_{ij})^2}$, where $T$ is the test set.

\begin{figure}[!t]
\begin{center}
\includegraphics[width=0.495\columnwidth]{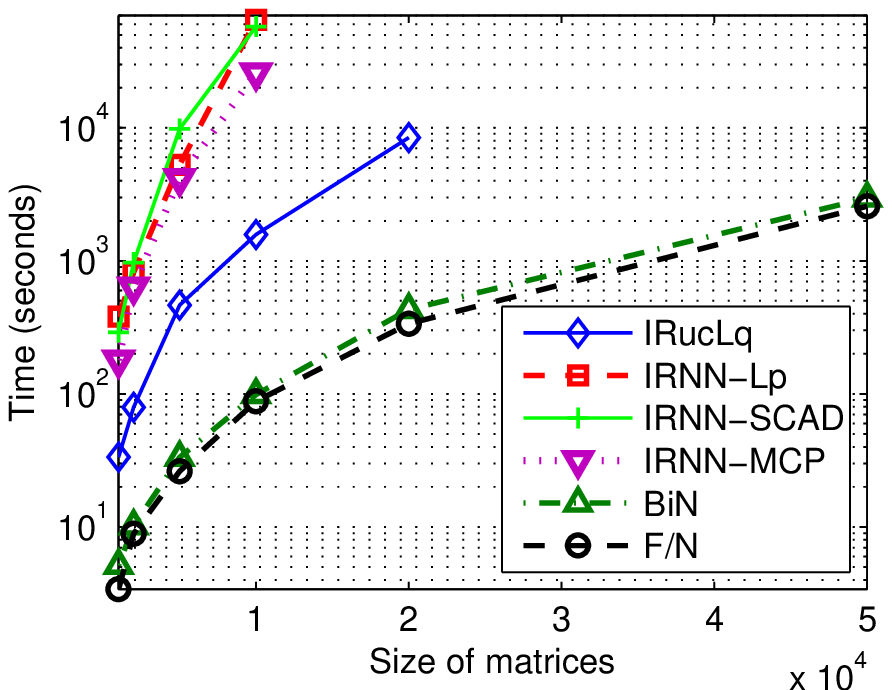}
\includegraphics[width=0.495\columnwidth]{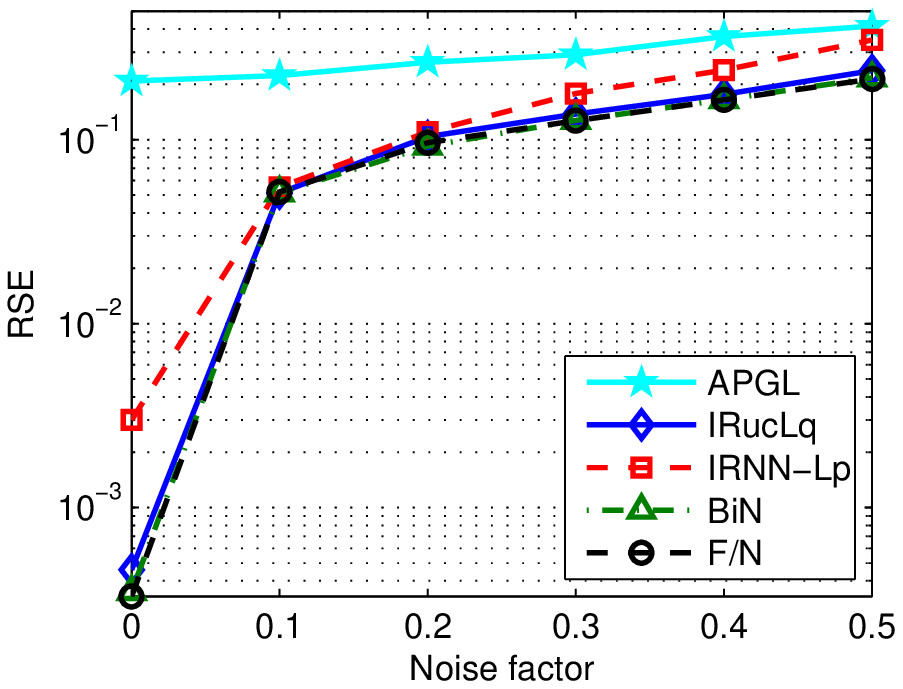}
\caption{The running time (seconds) and RSE results vs. sizes of matrices (left) and noise factors (right).}
\label{fig_sim2}
\end{center}
\end{figure}

\subsection{Synthetic Matrix Completion}
The synthetic matrices $Z\!\in\!\mathbb{R}^{m\times n}$ with rank $r$ are generated by the following procedure: the entries of both $U\!\in\!\mathbb{R}^{m\times r}$ and $V\!\in\!\mathbb{R}^{n\times r}$ are first generated as independent and identically distributed (i.i.d.) numbers, and then $Z\!=\!UV^{T}$ is assembled. Since all these algorithms have very similar recovery performance on noiseless matrices, we only conducted experiments on noisy matrices with different noise levels, i.e., $\mathcal{P}_{\Omega}(D)\!=\!\mathcal{P}_{\Omega}(Z\!+\!nf\!*\!E)$, where $nf$ denotes the noise factor. In other worlds, the observed subset is corrupted by i.i.d.\ standard Gaussian random noise as in~\citep{lu:lrm}. In addition, only 20\% or 30\% entries of $D$ are sampled uniformly at random as training data, i.e., sampling ratio (SR)$=$20\% or 30\%. The rank parameter $d$ of our algorithms is set to $\lfloor1.25r\rfloor$ as in~\citep{wen:nsor}.

The average RSE results of 100 independent runs on noisy random matrices are shown in Figure 1, which shows that if $p$ varies from 0.1 to 0.7, IRucLq and IRNN-Lp achieve similar recovery performance as IRNN-SCAD, IRNN-MCP and our algorithms; otherwise, IRucLq and IRNN-Lp usually perform much worse than the other four methods, especially $p=1$. We also report the running time of all the methods with 20\% SR as the size of noisy matrices increases, as shown in Figure 2. Moreover, we present the RSE results of those methods and APGL{\footnote{\url{http://www.math.nus.edu.sg/~mattohkc/}}}~\citep{toh:apg} (which is one of the nuclear norm solvers) with different noise factors. Figure 2 shows that our algorithms are significantly faster than the other methods, while the running time of IRucLq and IRNN increases dramatically when the size of matrices increases, and they could not yield experimental results within 48 hours when the size of matrices is $50,000\!\times\!50,000$. This further justifies that both our algorithms have very good scalability and can address large-scale problems. In addition, with only 20\% SR, all Schatten quasi-norm methods significantly outperform APGL in terms of RSE.

\begin{table}[!t]
\centering
\setlength{\tabcolsep}{2.5pt}
\scriptsize
\caption{Testing RMSE on MovieLens1M, MovieLens10M and Netflix.}
\label{tab_sim1}
\begin{tabular}{l|c|c|c}
\hline
{Datasets} & {MovieLens1M} &{MovieLens10M} & {Netflix}\\
\hline
{\% SR}  & 50\%\;/\;70\%\;/\;90\%     & 50\%\;/\;70\%\;/\;90\%       & 50\%\;/\;70\%\;/\;90\%\\
\hline
{APGL}       &1.2564/\,1.1431/\,0.9897	&1.1138/\,0.9455/\,0.8769	&1.0806/\,0.9885/\,0.9370\\
{LMaFit}     &0.9138/\,0.9019/\,0.8845	&0.8705/\,0.8496/\,0.8244	&0.9062/\,0.8923/\,0.8668\\
{IRucLq}     &0.9099/\,0.8918/\,0.8786	&---\;/\;---\;/\;---	&---\;/\;---\;/\;---\\
{IRNN}       &0.9418/\,0.9275/\,0.9032	&---\;/\;---\;/\;---	&---\;/\;---\;/\;---\\
{BiN}        &\textbf{0.8741}/\,0.8593/\,0.8485  &0.8274/\,0.8115/\,0.7989  &0.8650/\,0.8487/\,0.8413\\
{F/N}        &0.8764/\,\textbf{0.8562}/\,\textbf{0.8441} &\textbf{0.8158}/\,\textbf{0.8021}/\,\textbf{0.7921} &\textbf{0.8618}/\,\textbf{0.8459}/\,\textbf{0.8404}\\
\hline
\end{tabular}
\end{table}

\begin{figure*}[t]
\begin{center}
\subfigure[MovieLens1M]{\includegraphics[width=0.609\columnwidth]{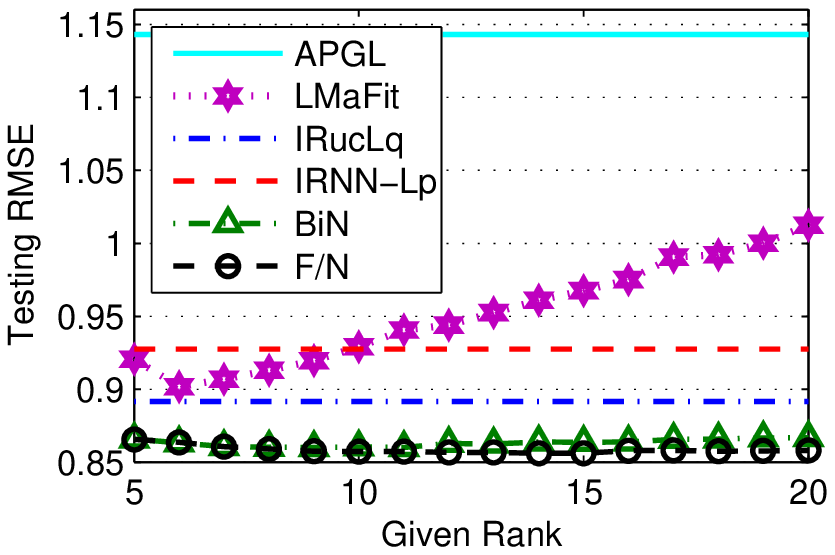}}\;\;\;\;
\subfigure[MovieLens10M]{\includegraphics[width=0.609\columnwidth]{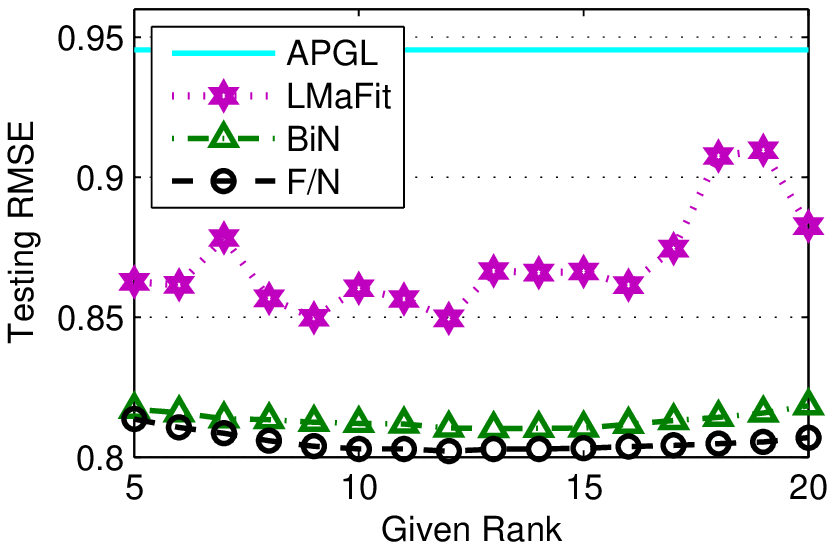}}\;\;\;\;
\subfigure[Netflix]{\includegraphics[width=0.609\columnwidth]{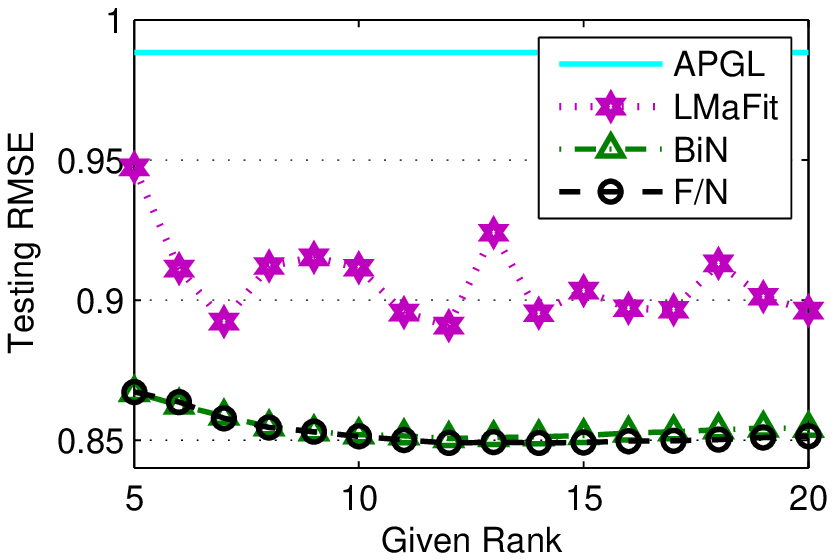}}
\vspace{-2mm}
\caption{The testing RMSE of LMaFit and our algorithms with ranks varying from 5 to 20 and 70\% SR.}
\label{fig_sim3}
\end{center}
\end{figure*}

\begin{figure}[!t]
\begin{center}
\includegraphics[width=0.469\columnwidth]{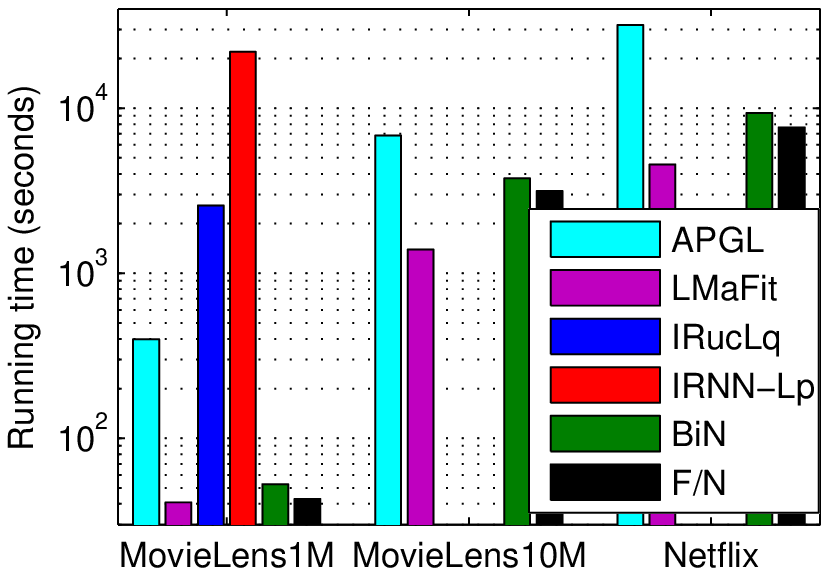}
\includegraphics[width=0.479\columnwidth]{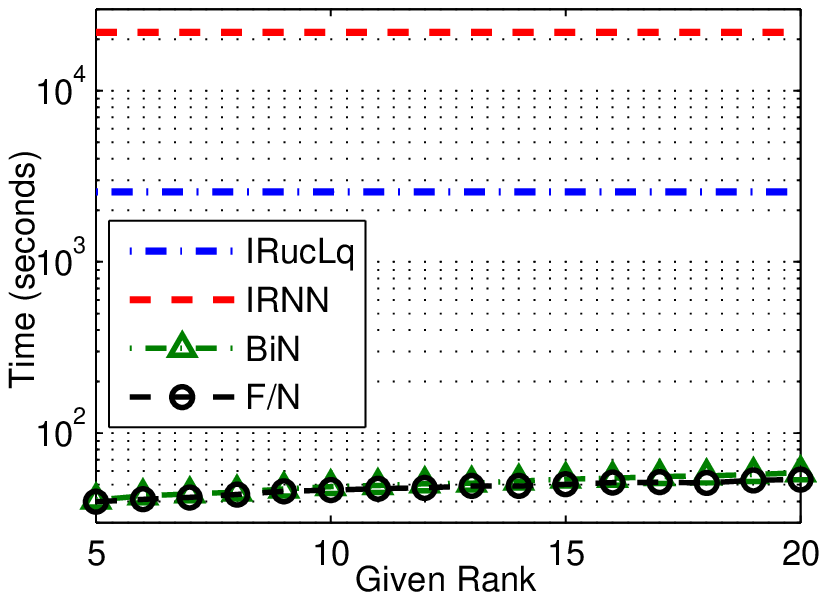}
\vspace{-2mm}
\caption{The running time (seconds) on three data sets (left, best viewed in colors) and MovieLens1M (right).}
\label{fig_sim4}
\end{center}
\end{figure}

\subsection{Collaborative Filtering}
We tested our algorithms on three real-world recommendation system data sets: the MovieLens1M, MovieLens10M\footnote{\url{http://www.grouplens.org/node/73}} and  Netflix datasets~\citep{kdd:cup}. We randomly chose 50\%, 70\% and 90\% as the training set and the remaining as the testing set, and the experimental results are reported over 10 independent runs. In addition to the methods used above, we also compared our algorithms with one of the fastest existing methods, LMaFit{\footnote{\url{http://lmafit.blogs.rice.edu/.}}}~\citep{wen:nsor}. The testing RMSE of all these methods on the three data sets is reported in Table 1, which shows that all those methods with non-convex penalty functions perform significantly better than the convex nuclear norm solver, APGL. In addition, our algorithms consistently outperform the other methods in terms of prediction accuracy. This further confirms that our two Schatten quasi-norm regularized models can provide a good estimation of a low-rank matrix. Moreover, we report the average testing RMSE and running time of our algorithms on these three data sets in Figures 3 and 4, where the rank varies from 5 to 20 and SR is set to 70\%. Note that IRucLq and IRNN-Lp could not run on the two larger data sets due to runtime exceptions. It is clear that our algorithms are much faster than AGPL, IRucLq and IRNN-Lp on all these data sets. They perform much more robust with respect to ranks than LMaFit, and are comparable in speed with it. This shows that our algorithms have very good scalability and are suitable for real-world applications.

\begin{figure}[!t]
\small
\centering
\subfigure[Original]{\includegraphics[width=0.24\linewidth]{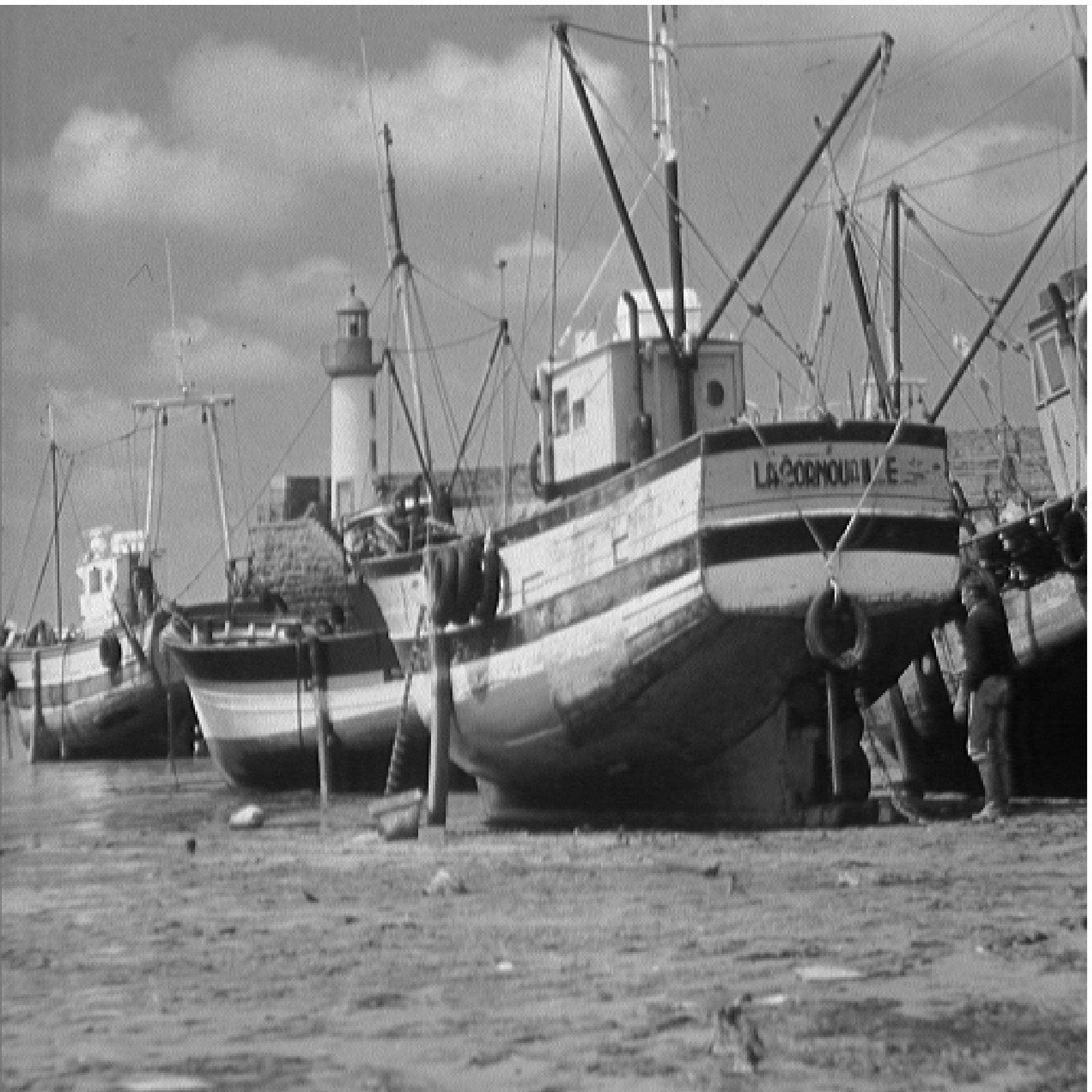}}
\subfigure[Input]{\includegraphics[width=0.24\linewidth]{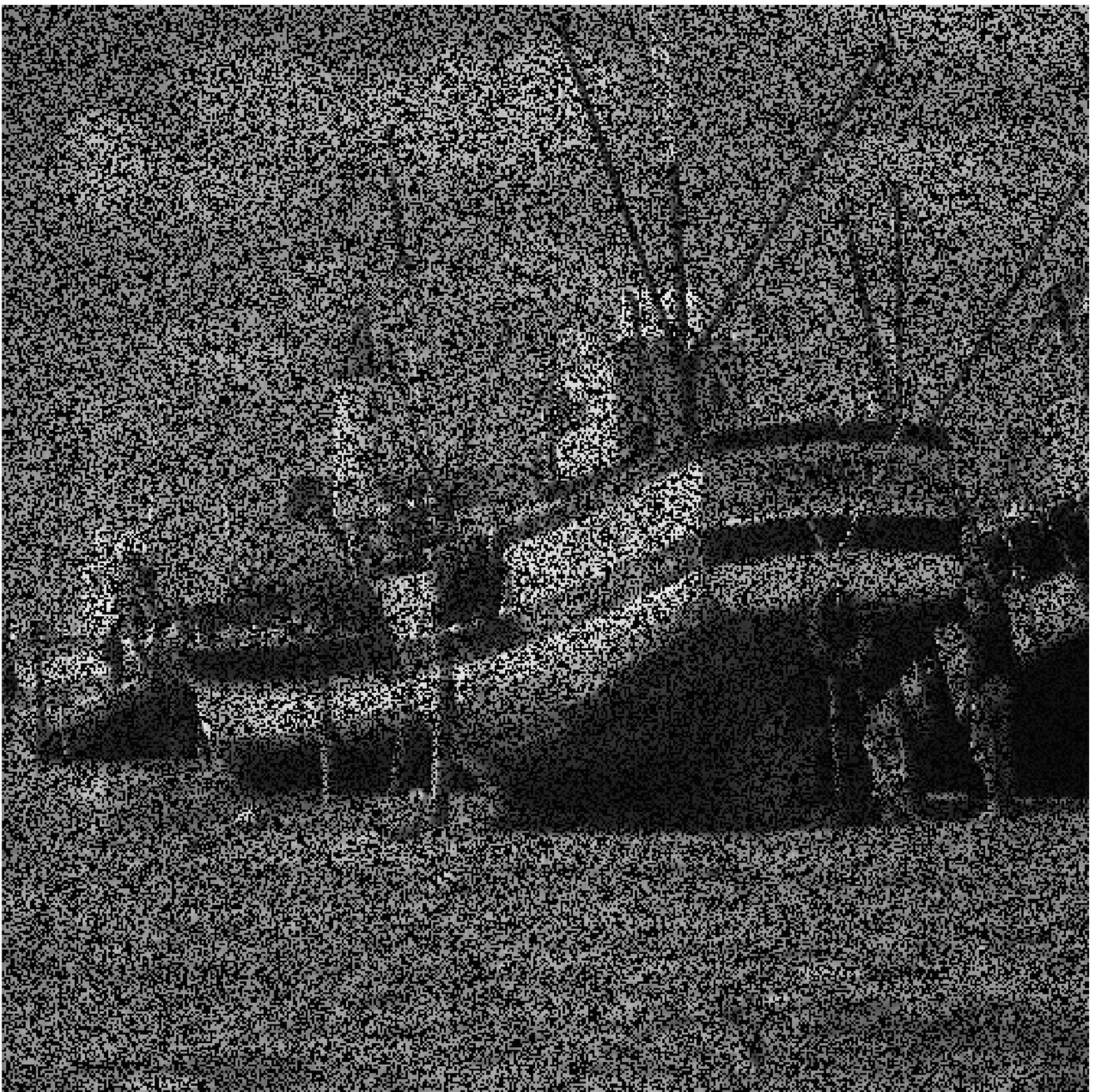}}
\subfigure[APGL]{\includegraphics[width=0.24\linewidth]{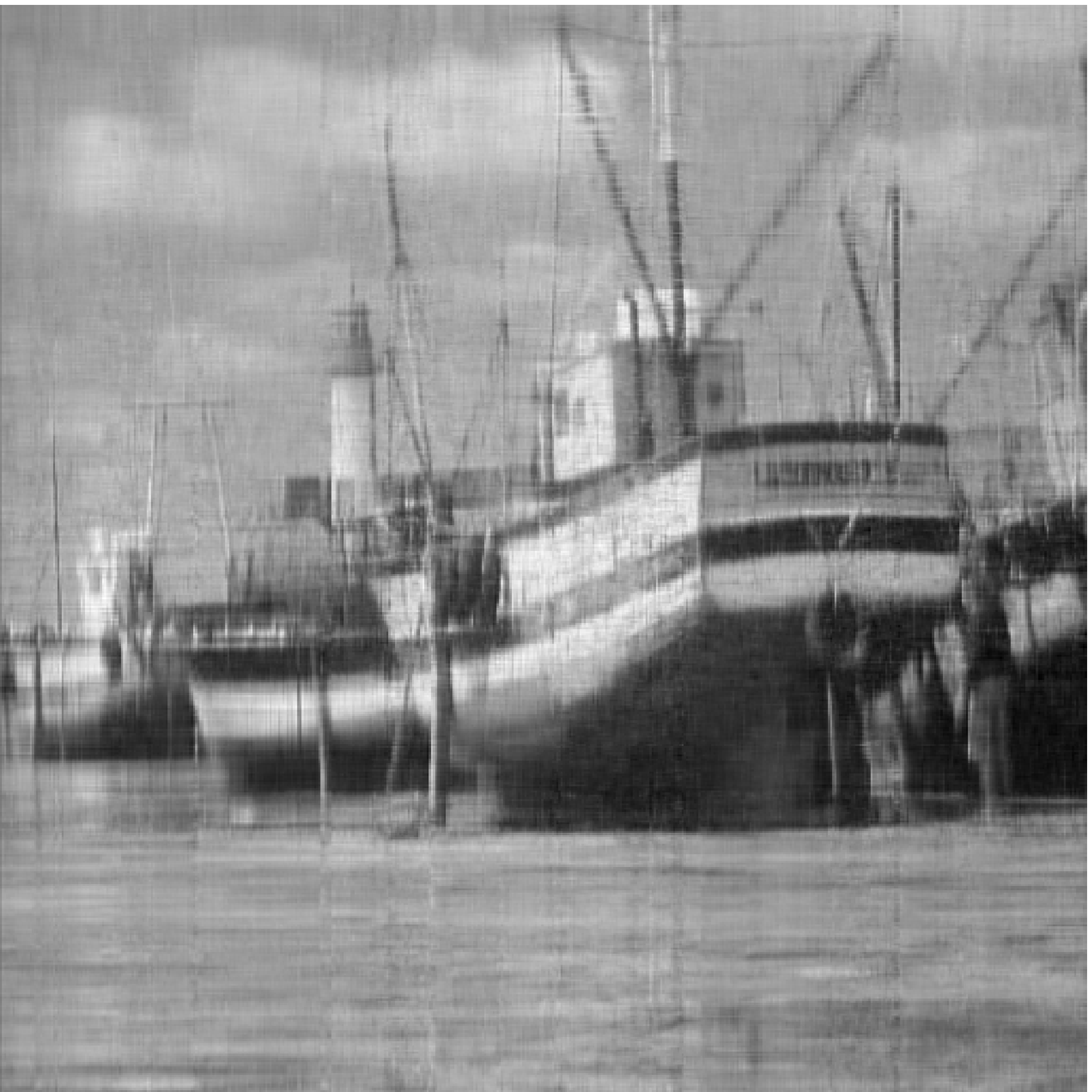}}
\subfigure[LMaFit]{\includegraphics[width=0.24\linewidth]{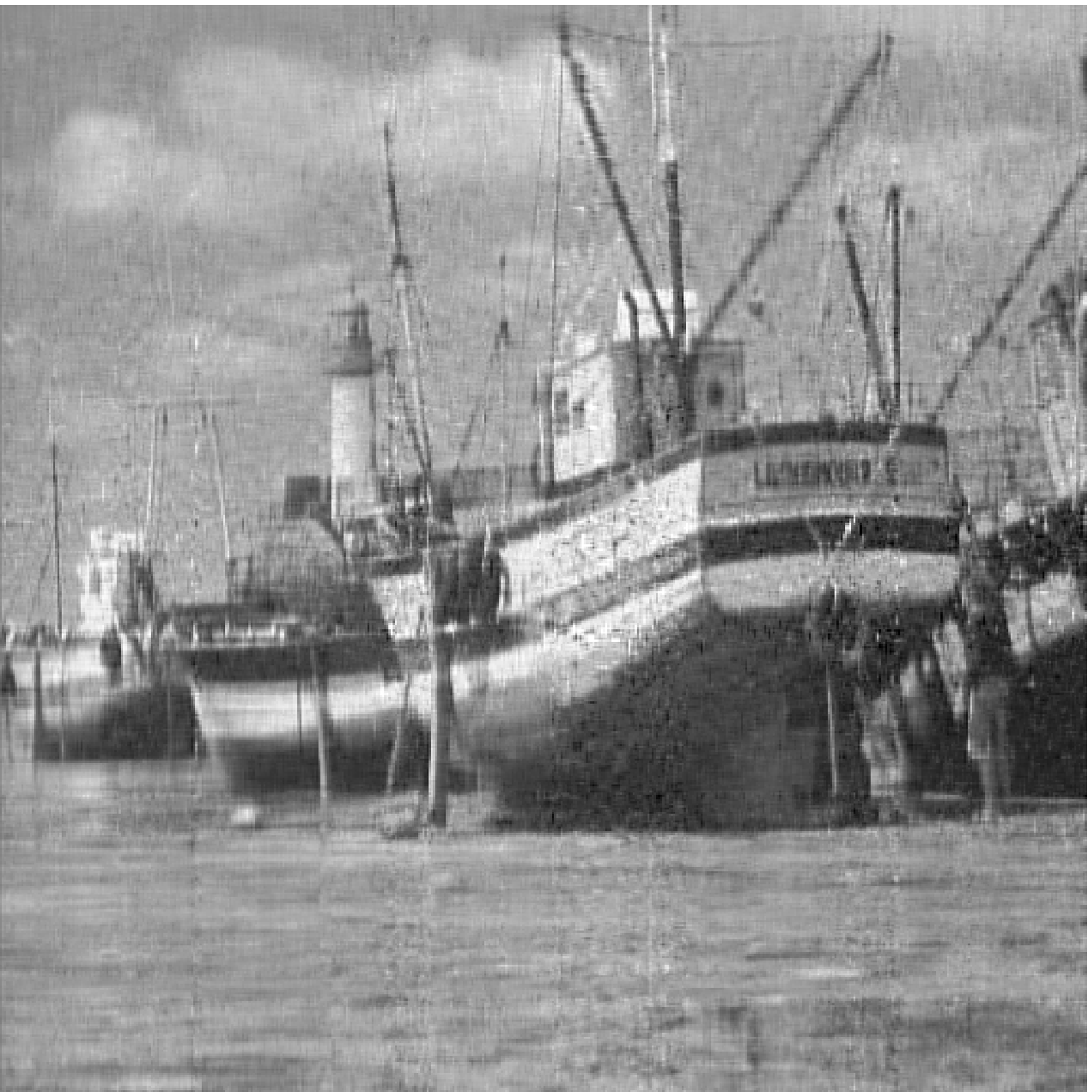}}
\subfigure[IRucLq]{\includegraphics[width=0.24\linewidth]{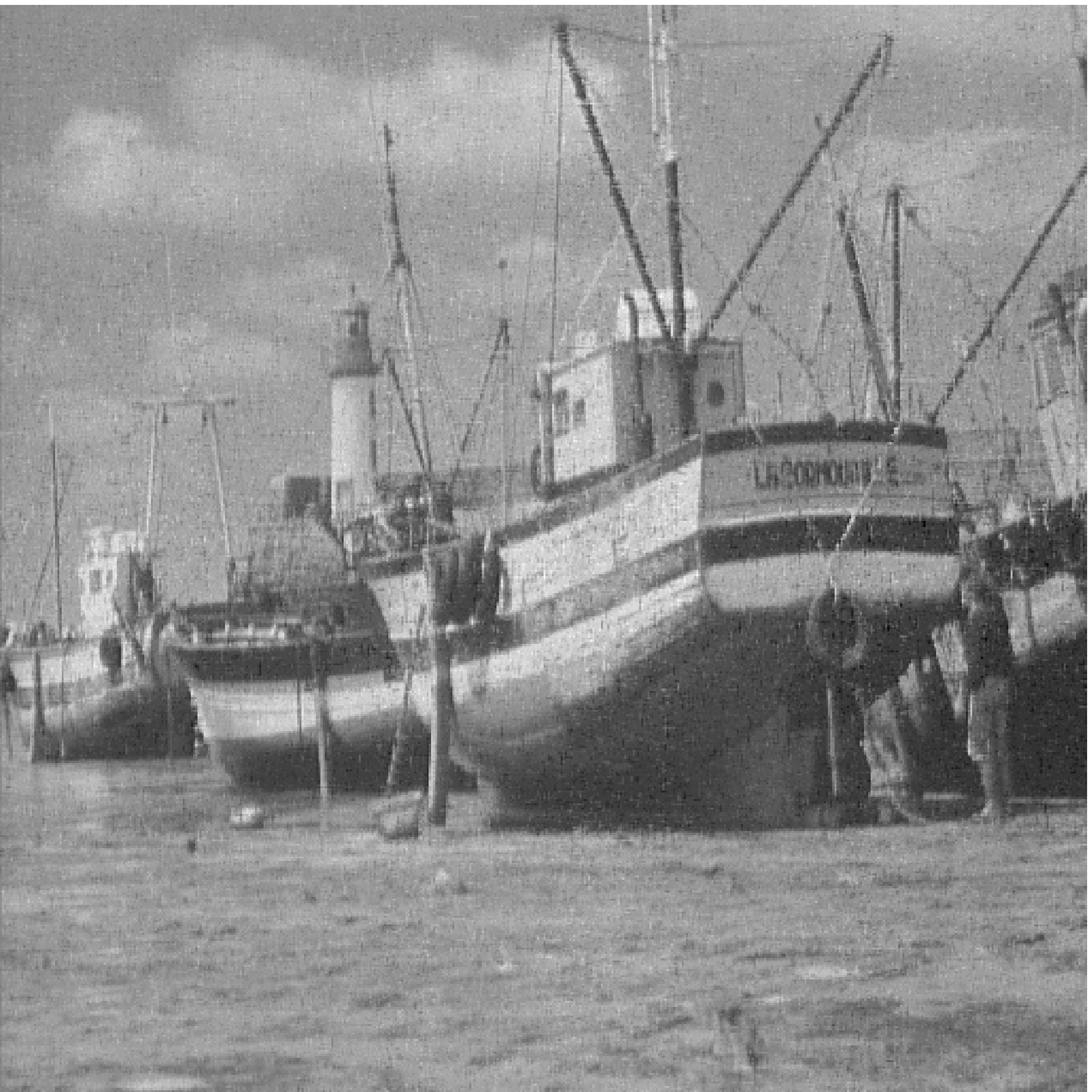}}
\subfigure[IRNN-Lp]{\includegraphics[width=0.24\linewidth]{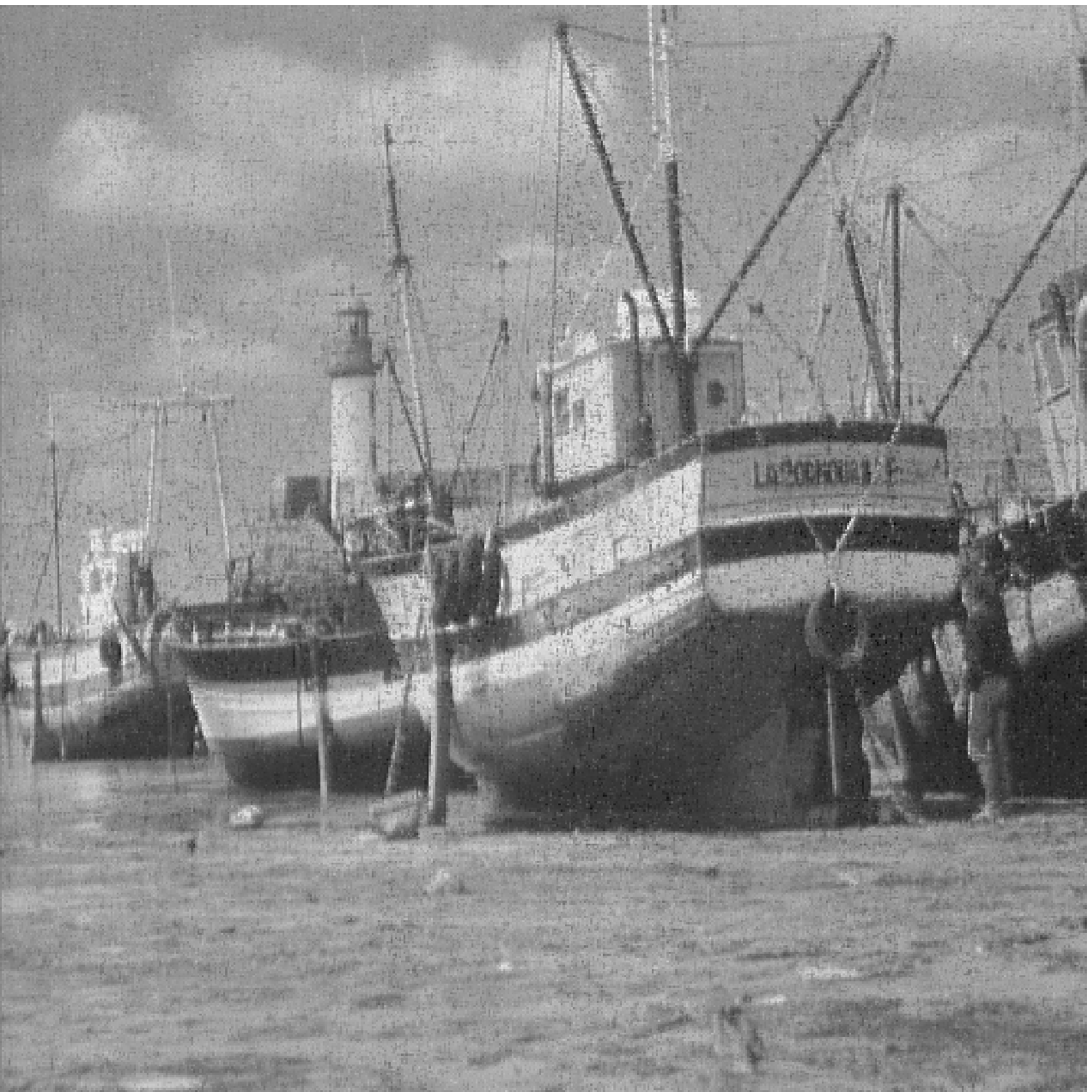}}
\subfigure[BiN]{\includegraphics[width=0.24\linewidth]{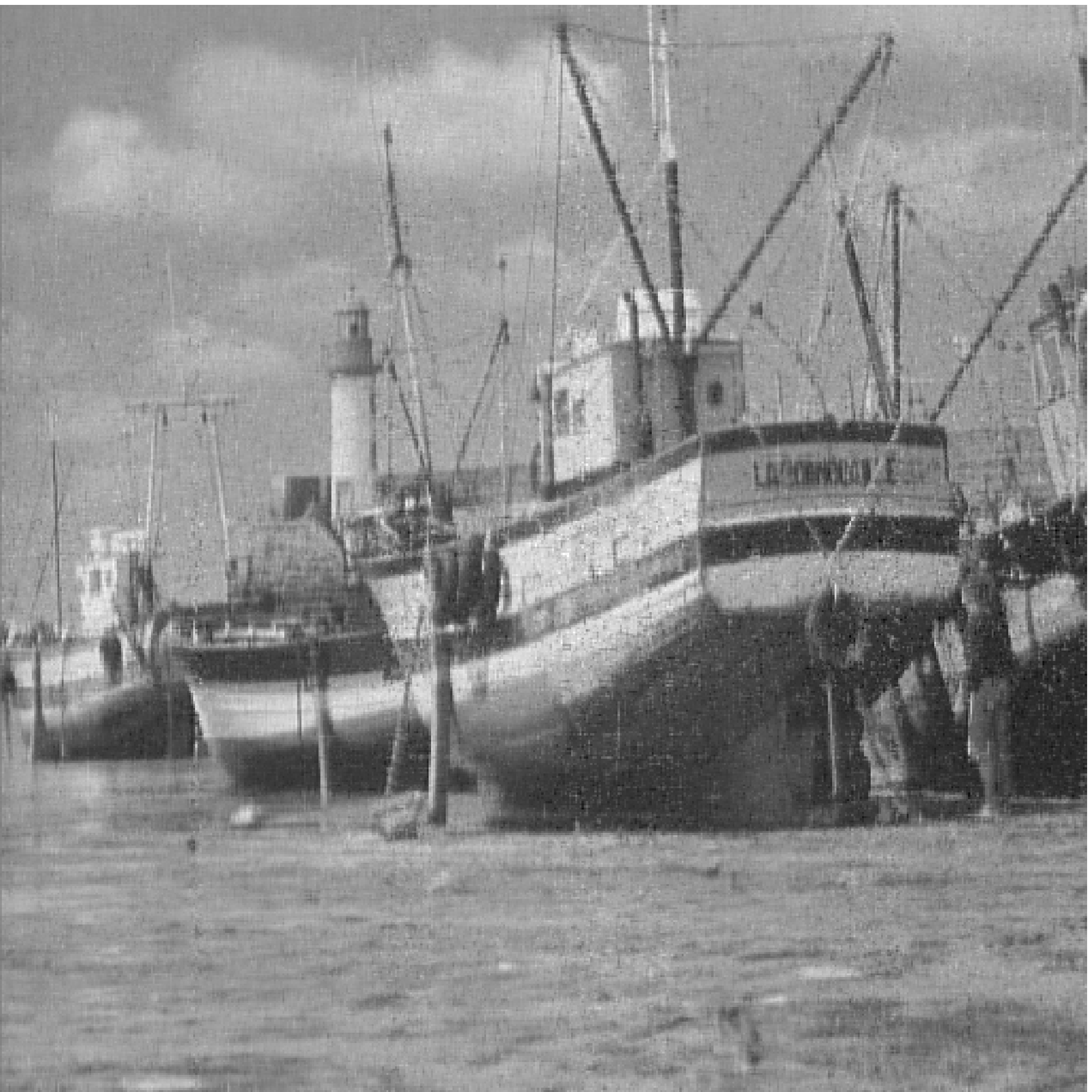}}
\subfigure[F/N]{\includegraphics[width=0.24\linewidth]{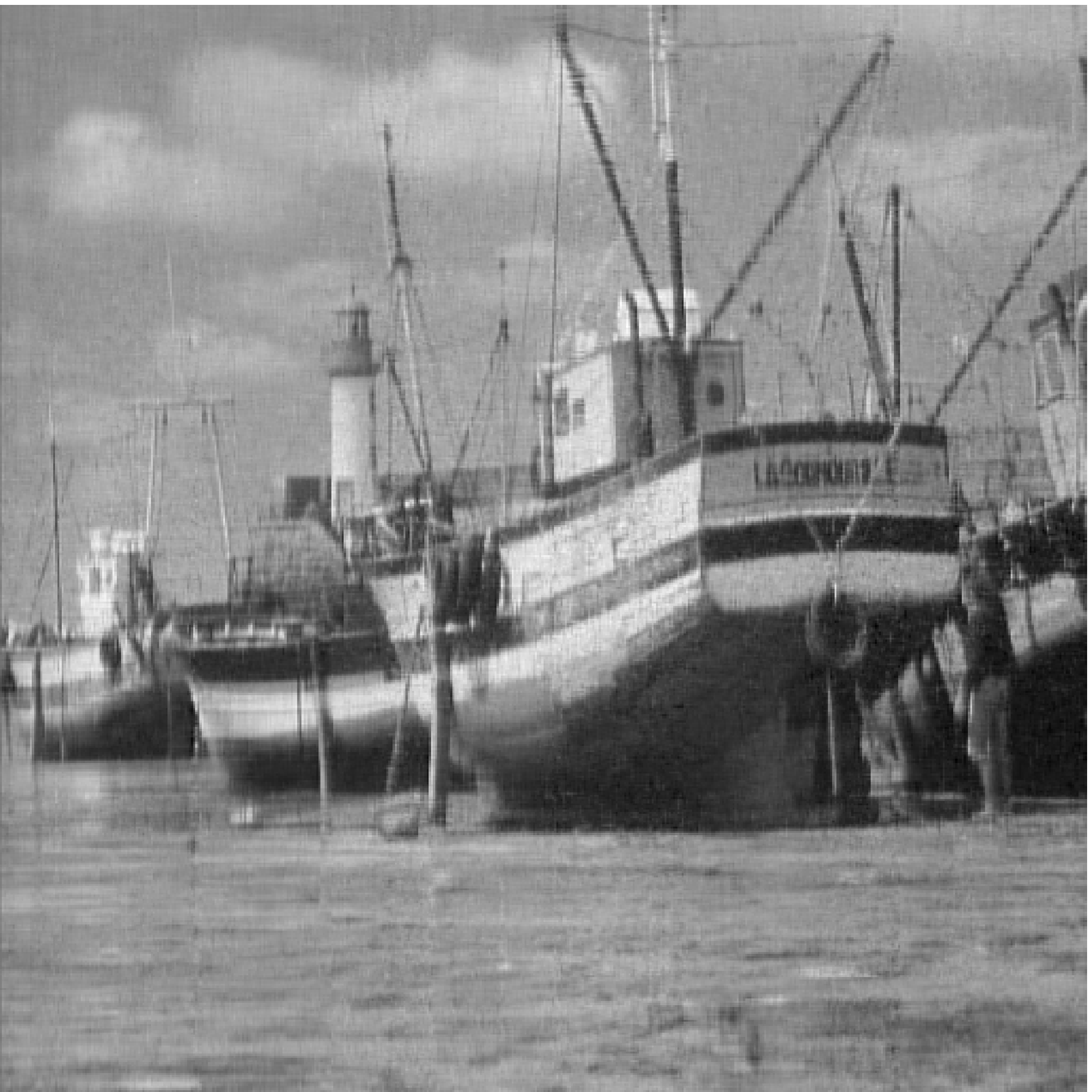}}
\vspace{-2mm}
\caption{Comparison of image recovery on the Boat image of size $512\!\times\!512$: (a) Original image; (b) Image with Gaussian noise; (c) APGL (PSNR: 24.93, Time: 15.47sec); (d) LMaFit (PSNR: 25.89, Time: 6.95sec); (e) IRucLq (PSNR: 26.36, Time: 805.81sec); (f) IRNN-Lp (PSNR: 26.21, Time: 943.28sec); (g) BiN (PSNR: 26.94, Time: 8.93sec); (h) F/N (PSNR: 27.62, Time: 10.80sec).}
\label{fig_sim5}
\end{figure}

\subsection{Image Recovery}
We also applied our algorithms to gray-scale image recovery on the Boat image of size $512\!\times\!512$, where 50\% of pixels in the input image were replaced by random Gaussian noise, as shown in Figure 5(b). In addition, we employed the well known peak signal-to-noise ratio (PSNR) to measure the recovery performance. The rank parameter of our algorithms and IRucLq was set to 100. Due to limited space, we only report the best results (PSNR and CPU time) of APGL, LMaFit, IRucLq and IRNN-Lp in Figure 5, which shows that our two algorithms achieve much better recovery performance than the other methods in terms of PSNR. And impressively, both our algorithms are significantly faster than the other methods except LMaFit and at least 70 times faster than IRucLq and IRNN-Lp.

\section{Conclusions}
In this paper we defined two tractable Schatten quasi-norms, i.e., the Frobenius/nuclear hybrid and bi-nuclear quasi-norms, and proved that they are in essence the Schatten-$2/3$ and $1/2$ quasi-norms, respectively. Then we designed two efficient proximal alternating linearized minimization algorithms to solve our Schatten quasi-norm minimization for matrix completion problems, and also proved that each bounded sequence generated by our algorithms globally converges to a critical point. In other words, our algorithms not only have better convergence properties than existing algorithms, e.g., IRucLq and IRNN, but also reduce the computational complexity from $O(mn^{2})$ to $O(mnd)$, with $d$ being the estimated rank ($d\!\ll\! m,n$). We also provided the recovery guarantee for our algorithms, which implies that they need only $O(md\log(m))$ observed entries to recover a low-rank matrix with high probability. Our experiments showed that our algorithms outperform the state-of-the-art methods in terms of both efficiency and effectiveness.

\section{Acknowledgements}
We thank the reviewers for their constructive comments. The authors are partially supported by the SHIAE fund 8115048 and the Hong Kong GRF 2150851.

\fontsize{9.2pt}{10.2pt}
\selectfont
\bibliography{aaai16_1718}
\bibliographystyle{aaai}

\end{document}


%
\title{Supplementary Materials for ``Scalable Algorithms for Tractable Schatten Quasi-Norm Minimization"}

\author{Fanhua Shang, Yuanyuan Liu, James Cheng\\
Department of Computer Science and Engineering, The Chinese University of Hong Kong\\
\{fhshang, yyliu, jcheng\}@cse.cuhk.edu.hk
}
\maketitle

\setcounter{equation}{9}
\setcounter{lemma}{1}
\setcounter{definition}{3}

In this supplementary material, we give the detailed proofs of some properties and theorems.

\section{More Notations}
$\mathbb{R}^{n}$ denotes the $n$-dimensional Euclidean space, and the set of all $m\times n$ matrices with real entries is denoted by $\mathbb{R}^{m\times n}$. Given two matrices, $X$ and $Y\in\mathbb{R}^{m\times n}$, the inner product is defined by $\langle X, Y\rangle:=\textrm{Tr}(X^{T}Y)$, where $\textrm{Tr}(\cdot)$ denotes the trace of a matrix. $\|X\|_{2}$ is the spectral norm and is equal to the maximum singular value of $X$. $I$ denotes an identity matrix.

For any vector $x\in\mathbb{R}^{n}$, its $l_{p}$ quasi-norm for $0\!<\!p\!<\!1$ is defined as
\begin{displaymath}
\|x\|_{p}=\left(\sum_{i}|x_{i}|^{p}\right)^{1/p}.
\end{displaymath}
In addition, the $l_{1}$-norm and the $l_{2}$-norm of $x$ are $\|x\|_{1}=\sum_{i}|x_{i}|$ and $\|x\|_{2}=\sqrt{\sum_{i}x^{2}_{i}}$, respectively.

For any matrix $X\in\mathbb{R}^{m\times n}$, we assume the singular values of $X$ are ordered as $\sigma_{1}(X)\geq \sigma_{2}(X)\geq\cdots\geq\sigma_{r}(X)>\sigma_{r+1}(X)=\cdots=\sigma_{\min(m,n)}(X)=0$, where $r=\textrm{rank}(X)$. By writing $X=U\Sigma V^{T}$ in its standard singular value decomposition (SVD), we can extend $X=U\Sigma V^{T}$ to the following definitions.

The Schatten-$p$ quasi-norm ($0<p<1$) of a matrix $X\in\mathbb{R}^{m\times n}$ is defined as follows:
\begin{displaymath}
\|X\|_{S_{p}}=\left(\sum^{\min(m,n)}_{i=1}\left(\sigma_{i}(X)\right)^{p}\right)^{1/p}.
\end{displaymath}

The nuclear norm (also called the trace norm or the Schatten-1 norm) of $X$ is defined as
\begin{displaymath}
\|X\|_{*}=\sum^{\min(m,n)}_{i=1}\sigma_{i}(X).
\end{displaymath}

The Frobenius norm (also called the Schatten-2 norm) of $X$ is defined as
\begin{displaymath}
\|X\|_{F}=\sqrt{\textrm{Tr}\left(X^{T}X\right)}=\sqrt{\sum^{\min(m,n)}_{i=1}\sigma^{2}_{i}(X)}.
\end{displaymath}

\section{Proofs of Theorem 1 and Property 1:}
In the following, we will first prove that the Frobenius/nuclear hybrid norm $\|\!\cdot\!\|_{\textup{F/N}}$ is a quasi-norm.
\begin{proof}
By the definition of the F/N norm, for any $a$, $a_{1}$, $a_{2}\in \mathbb{R}$ and $a=a_{1}a_{2}$, we have
\begin{displaymath}
\begin{split}
\|aX\|_{\textup{F/N}}&=\min_{aX=(a_{1}U)(a_{2}V^{T})}\|a_{1}U\|_{*}\|a_{2}V\|_{F}\\
&=\min_{X=UV^{T}}(|a_{1}|\,\|U\|_{*})\,(|a_{2}|\,\|V\|_{F})\\
&=|a|\min_{X=UV^{T}}\|U\|_{*}\|V\|_{F}\\
&=|a|\,\|X\|_{\textup{F/N}}.
\end{split}
\end{displaymath}

Next we will prove that $\|X+Y\|_{\textup{F/N}}\leq \beta\left(\|X\|_{\textup{F/N}}+\|Y\|_{\textup{F/N}}\right)$, where $\beta\geq 1$. By Lemma 1, i.e., $\|X\|_{*}=\min_{X=UV^{T}}\|U\|_{F}\|V\|_{F}$, there must exist both matrices $\widehat{U}$ and $\widehat{V}$ such that $\|X\|_{*}=\|\widehat{U}\|_{F}\|\widehat{V}\|_{F}$ with the constraint $X=\widehat{U}\widehat{V}$. According to the definition of the F/N-norm and the fact that $\|X\|_{*}\leq\sqrt{\textrm{rank}(X)}\|X\|_{F}$, we have
\begin{displaymath}
\begin{split}
\|X\|_{\textup{F/N}}&=\min_{X=UV^{T}}\|U\|_{*}\|V\|_{F}\\
&\leq \|\widehat{U}\|_{*}\|\widehat{V}\|_{F}\\
&\leq \sqrt{\textrm{rank}(U)}\|\widehat{U}\|_{F}\|\widehat{V}\|_{F}\\
&\leq \sqrt{\textrm{rank}(U)}\|X\|_{*}\\
&=\alpha\|X\|_{*},
\end{split}
\end{displaymath}
where $\alpha=\sqrt{\textrm{rank}(U)}$. If $X\neq0$, we can know that $\alpha\geq 1$. On the other hand, we also have
\begin{displaymath}
\|X\|_{*}\leq \|X\|_{\textup{F/N}}.
\end{displaymath}
By the above properties, there exists a constant $\beta\geq 1$ such that the following holds for all $X,Y\in\mathbb{R}^{m\times n}$
\begin{displaymath}
\begin{split}
\|X+Y\|_{\textup{F/N}}&\leq \beta\|X+Y\|_{*}\\
&\leq \beta(\|X\|_{*}+\|Y\|_{*})\\
&\leq \beta(\|X\|_{\textup{F/N}}+\|Y\|_{\textup{F/N}}).
\end{split}
\end{displaymath}

Furthermore, $\forall X\in\mathbb{R}^{m\times n}$ and $X=UV^{T}$, we have
\begin{displaymath}
\|X\|_{\textup{F/N}}=\min_{X=UV^{T}}\|U\|_{*}\|V\|_{F}\geq0.
\end{displaymath}
In addition, if $\|X\|_{\textup{F/N}}=0$, we have $\|X\|_{*}\leq\|X\|_{\textup{F/N}}=0$, i.e., $\|X\|_{*}=0$.
Hence, $X=0$. In short, the F/N-norm $\|\cdot\|_{\textup{F/N}}$ is a quasi-norm.
\end{proof}

\onecolumn
Before giving the complete proofs for Theorem 1 and Property 1, we first present and prove the following important lemma.

\begin{lemma}
Suppose that $Z\in\mathbb{R}^{m\times n}$ is a matrix of rank $r\leq\min(m,\,n)$, and we denote its SVD by $Z=L\Sigma_{Z}R^{T}$, where $L\in \mathbb{R}^{m\times r}$, $R\in \mathbb{R}^{n\times r}$ and $\Sigma_{Z}\in \mathbb{R}^{r\times r}$. For any matrix $A\in\mathbb{R}^{r\times r}$ satisfying $AA^{T}=A^{T}A=I_{r\times r}$, and the given $p\;(0<p<1)$, then $(A\Sigma_{Z}A^{T})_{k,k}\geq 0$ for all $k=1,\ldots,r$, and
\begin{equation*}
\textup{Tr}^{p}(A\Sigma_{Z}A^{T})\geq \textup{Tr}^{p}(\Sigma_{Z})=\|Z\|^{p}_{S_{p}},
\end{equation*}
where $\textup{Tr}^{p}(B)=\sum_{i}B^{p}_{ii}$.
\end{lemma}
\begin{proof}
For any $k\in\{1,\ldots,r\}$, we have $(A\Sigma_{Z}A^{T})_{k,k}=\sum_{i}a^{2}_{ki}\sigma_{i}\geq 0$, where $\sigma_{i}\geq 0$ is the $i$-th singular value of $Z$. Then
\begin{equation}\label{Seq1}
\textup{Tr}^{p}(A\Sigma_{Z}A^{T})=\sum_{k}\left(\sum_{i}a^{2}_{ki}\sigma_{i}\right)^{p}.
\end{equation}
Since $\psi(x)=x^{p}\;(0<p<1)$ is a concave function on $\mathbb{R}^{+}$, and by the Jensen's inequality~\citep{mitrinovic:ai} and $\sum_{i}a^{2}_{ki}=1$ for any $k\in\{1,\ldots,r\}$, we have
\begin{equation*}
\left(\sum_{i}a^{2}_{ki}\sigma_{i}\right)^{p}\geq \sum_{i}a^{2}_{ki}\sigma^{p}_{i}.
\end{equation*}
According to the above inequality and $\sum_{k}a^{2}_{ki}=1$ for any $i\in\{1,\ldots,r\}$, \eqref{Seq1} can be rewritten as
\begin{equation*}
\begin{split}
\textup{Tr}^{p}(A\Sigma_{Z}A^{T}) &= \sum_{k}\left(\sum_{i}a^{2}_{ki}\sigma_{i}\right)^{p} \\
&\geq \sum_{k}\sum_{i}a^{2}_{ki}\sigma^{p}_{i}\\
&=\sum_{i}\sigma^{p}_{i}\\
&=\textup{Tr}^{p}(\Sigma_{Z})=\|Z\|^{p}_{S_{p}}.
\end{split}
\end{equation*}
This completes the proof.
\end{proof}

\emph{\textbf{Proof of Theorem 1:}}
\begin{proof}
Assume that $U\!=\!L_{U}\Sigma_{U}R^{T}_{U}$ and $V\!=\!L_{V}\Sigma_{V}R^{T}_{V}$ are the thin SVDs of $U$ and $V$, respectively, where $L_{U}\in \mathbb{R}^{m\times d}$, $L_{V}\in \mathbb{R}^{n\times d}$, and $R_{U},\Sigma_{U},R_{V},\Sigma_{V}\in \mathbb{R}^{d\times d}$. Let $X\!=\!L_{X}\Sigma_{X}R^{T}_{X}$, where the columns of $L_{X}\in \mathbb{R}^{m\times d}$ and $R_{X}\in \mathbb{R}^{n\times d}$ are the left and right singular vectors associated with the top $d$ singular values of $X$ with rank at most $r$ $(r\!\leq\! d)$, and $\Sigma_{X}=\textup{diag}([\sigma_{1}(X),\!\cdots\!,\sigma_{r}(X),0,\!\cdots\!,0])\in\!\mathbb{R}^{d\times d}$.

By $X=UV^{T}$, i.e., $L_{X}\Sigma_{X}R^{T}_{X}=L_{U}\Sigma_{U}R^{T}_{U}R_{V}\Sigma_{V}L^{T}_{V}$, then $\exists O_{1},\widehat{O}_{1}\in \mathbb{R}^{d\times d}$ satisfy $L_{X}=L_{U}O_{1}$ and $L_{U}=L_{X}\widehat{O}_{1}$, i.e., $O_{1}=L^{T}_{U}L_{X}$ and $\widehat{O}_{1}=L^{T}_{X}L_{U}$. Thus, $O_{1}=\widehat{O}^{T}_{1}$. Since $L_{X}=L_{U}O_{1}=L_{X}\widehat{O}_{1}O_{1}$, we have $\widehat{O}_{1}O_{1}=O^{T}_{1}O_{1}=I_{d\times d}$. Similarly, we have $O_{1}\widehat{O}_{1}=O_{1}O^{T}_{1}=I_{d\times d}$. In addition, $\exists O_{2}\in \mathbb{R}^{d\times d}$ satisfies $R_{X}=L_{V}O_{2}$ with $O_{2}O^{T}_{2}=O^{T}_{2}O_{2}=I_{d\times d}$. Let $O_{3}=O_{2}O^{T}_{1}\in \mathbb{R}^{d\times d}$, then we have $O_{3}O^{T}_{3}=O^{T}_{3}O_{3}=I_{d\times d}$, i.e., $\sum_{i}(O_{3})^{2}_{ij}=\sum_{j}(O_{3})^{2}_{ij}=1$ for $\forall i,j\in \{1,2,\ldots,d\}$, where $a_{i,j}$ denotes the element of the matrix $A$ in the $i$-th row and the $j$-th column. Furthermore, let $O_{4}=R^{T}_{U}R_{V}$, we have $\sum_{i}(O_{4})^{2}_{ij}\leq 1$ and $\sum_{j}(O_{4})^{2}_{ij}\leq 1$ for $\forall i,j\in \{1,2,\ldots,d\}$.

By the above analysis, then we have $O_{2}\Sigma_{X}O^{T}_{2}=O_{2}O^{T}_{1}\Sigma_{U}O_{4}\Sigma_{V}=O_{3}\Sigma_{U}O_{4}\Sigma_{V}$. Let $\tau_{i}$ and $\varrho_{j}$ denote the $i$-th and the $j$-th diagonal elements of $\Sigma_{U}$ and $\Sigma_{V}$, respectively. By Lemma 2, and $p=2/3$, we have
\begin{equation*}
\begin{split}
\|X\|_{S_{2/3}}&\leq\left(\textup{Tr}^{\frac{2}{3}}(O_{2}\Sigma_{X}O^{T}_{2})\right)^{\frac{3}{2}}=\left(\textup{Tr}^{\frac{2}{3}}(O_{2}O^{T}_{1}\Sigma_{U}O_{4}\Sigma_{V})\right)^{\frac{3}{2}}=\left(\textup{Tr}^{\frac{2}{3}}(O_{3}\Sigma_{U}O_{4}\Sigma_{V})\right)^{\frac{3}{2}}\\
&=\left(
\sum^{d}_{i=1}\left[\sum^{d}_{j=1}\tau_{j}(O_{3})_{ij}(O_{4})_{ji}\varrho_{i}\right ]^{\frac{2}{3}}\right )^{\frac{3}{2}}\\
&=\left(
\sum^{d}_{i=1}\varrho^{\frac{2}{3}}_{i}\left(\sum^{d}_{j=1}\tau_{j}(O_{3})_{ij}(O_{4})_{ji}\right)^{\frac{2}{3}}\right )^{\frac{3}{2}}\\
&^{a}\!\!\!\leq \left(\left[
\sum^{d}_{i=1} (\varrho^{\frac{2}{3}}_{i})^{3}\right]^{\frac{1}{3}}\left[
\sum^{d}_{i=1} \left(\sum^{d}_{j=1}\tau_{j}(O_{3})_{ij}(O_{4})_{ji}\right)^{\frac{2}{3}\times \frac{3}{2}}\right]^{\frac{2}{3}}\right)^{\frac{3}{2}}\\
&=\sqrt{\sum^{d}_{i=1}\varrho^{2}_{i}}\sum^{d}_{i=1}\sum^{d}_{j=1}\tau_{j}(O_{3})_{ij}(O_{4})_{ji}\\
&^{b}\!\!\!\leq\sqrt{\sum^{d}_{i=1}\varrho^{2}_{i}}\sum^{d}_{i=1}\sum^{d}_{j=1}\tau_{j}\frac{(O_{3})^{2}_{ij}+(O_{4})^{2}_{ji}}{2}\\
&^{c}\!\!\!\leq\sqrt{\sum^{d}_{i=1}\varrho^{2}_{i}}\sum^{d}_{j=1}\tau_{j}\\
&= \|U\|_{*}\|V\|_{F}=\sqrt{\|U\|_{*}}\sqrt{\|U\|_{*}}\|V\|_{F}\\
&^{d}\!\!\!\leq \left( \frac{\sqrt{\|U\|_{*}}+\sqrt{\|U\|_{*}}+\|V\|_{F}}{3}\right)^{3}\\
&= \left( \frac{2\sqrt{\|U\|_{*}}+\sqrt{\|V\|^{2}_{F}}}{3}\right)^{3}\\
&^{e}\!\!\!\leq \left(\frac{2\|U\|_{*}+\|V\|^{2}_{F}}{3}\right)^{\frac{3}{2}}
\end{split}
\end{equation*}
where the inequality $^{a}\!\!\!\leq$ holds due to the H\"{o}lder's inequality~\citep{mitrinovic:ai}, i.e., $\sum^{n}_{k=1}|x_{k}y_{k}|\leq(\sum^{n}_{k=1}|x_{k}|^{p})^{1/p}(\sum^{n}_{k=1}|y_{k}|^{q})^{1/q}$ with $1/p+1/q=1$, and here we set $p\!=\!3$ and $q\!=\!3/2$ in the inequality $^{a}\!\!\!\leq$\,; the inequality $^{b}\!\!\!\leq$ follows from the basic inequality $xy\leq \frac{x^{2}+y^{2}}{2}$ for any real numbers $x$ and $y$; the inequality $^{c}\!\!\!\leq$ relies on the facts $\sum_{i}(O_{3})^{2}_{ij}=1$ and $\sum_{i}(O_{4})^{2}_{ji}\leq 1$; the inequality $^{d}\!\!\!\leq$ holds due to the fact $\sqrt[3]{x_{1}x_{2}x_{3}}\leq (|x_{1}|+|x_{2}|+|x_{3}|)/3$ and the inequality $^{e}\!\!\!\leq$ holds due to the Jensen's inequality for the concave function $f(x)=x^{1/2}$.
Thus, for any matrices $U\in\mathbb{R}^{m\times d}$ and $V\in\mathbb{R}^{n\times d}$ satisfying $X=UV^{T}$, we have
\begin{equation*}
\|X\|_{S_{2/3}}\leq\|U\|_{*}\|V\|_{F}\leq\left(\frac{2\|U\|_{*}+\|V\|^{2}_{F}}{3}\right)^{\frac{3}{2}}.
\end{equation*}

On the other hand, let $U_{\star}=L_{X}\Sigma^{2/3}_{X}$ and $V_{\star}=R_{X}\Sigma^{1/3}_{X}$, where $\Sigma^{p}$ is entry-wise power to $p$, then we have $X=U_{\star}V^{T}_{\star}$ and
\begin{equation*}
\|X\|_{S_{2/3}}=\left(\textup{Tr}^{2/3}(\Sigma_{X})\right)^{\frac{3}{2}}=\|U_{\star}\|_{*}\|V_{\star}\|_{F}=\left(\frac{2\|U_{\star}\|_{*}+\|V_{\star}\|^{2}_{F}}{3}\right)^{\frac{3}{2}}.
\end{equation*}
Therefore, under the constraint $X=UV^{T}$, we have
\begin{displaymath}
\begin{split}
\|X\|_{S_{2/3}}=\min_{X=UV^{T}}\|U\|_{*}\|V\|_{F}=\min_{X=UV^{T}}\left(\frac{2\|U\|_{*}+\|V\|^{2}_{F}}{3}\right)^{\frac{3}{2}}=\|X\|_{\textup{F/N}}.
\end{split}
\end{displaymath}
This completes the proof.
\end{proof}
~\\

\section{Proof of Theorem 2:}

In the following, we will first prove that the bi-nuclear norm $\|\!\cdot\!\|_{\textup{BiN}}$ is a quasi-norm.

\begin{proof}
By the definition of the bi-nuclear norm, for any $a$, $a_{1}$, $a_{2}\in \mathbb{R}$ and $a=a_{1}a_{2}$, we have
\begin{displaymath}
\begin{split}
\|aX\|_{\textup{BiN}}&=\min_{aX=(a_{1}U)(a_{2}V^{T})}\|a_{1}U\|_{*}\|a_{2}V\|_{*}\\
&=\min_{X=UV^{T}}|a|\,\|U\|_{*}\|V\|_{*}\\
&=|a|\min_{X=UV^{T}}\|U\|_{*}\|V\|_{*}\\
&=|a|\,\|X\|_{\textup{BiN}}.
\end{split}
\end{displaymath}

Since $\|X\|_{*}=\min_{X=UV^{T}}\|U\|_{F}\|V\|_{F}$, and by Lemma 6 in~\citep{mazumder:sr}, there exist both matrices $\widehat{U}=U_{X}\Sigma^{1/2}_{X}$ and $\widehat{V}=V_{X}\Sigma^{1/2}_{X}$ such that $\|X\|_{*}=\|\widehat{U}\|_{F}\|\widehat{V}\|_{F}$ with the SVD of $X$, i.e., $X=U_{X}\Sigma_{X}V^{T}_{X}$. By the fact that $\|X\|_{*}\leq\sqrt{\textrm{rank}(X)}\|X\|_{F}$, we have
\begin{displaymath}
\begin{split}
\|X\|_{\textup{BiN}}&=\min_{X=UV^{T}}\|U\|_{*}\|V\|_{*}\\
&\leq \|\widehat{U}\|_{*}\|\widehat{V}\|_{*}\\
&\leq \sqrt{\textrm{rank}(X)}\sqrt{\textrm{rank}(X)}\|\widehat{U}\|_{F}\|\widehat{V}\|_{F}\\
&\leq \textrm{rank}(X)\|X\|_{*}.
\end{split}
\end{displaymath}
If $X\neq0$, then $\textrm{rank}(X)\geq 1$. On the other hand, we also have
\begin{displaymath}
\|X\|_{*}\leq \|X\|_{\textup{BiN}}.
\end{displaymath}
By the above properties, there exists a constant $\alpha\geq 1$ such that the following holds for all $X,Y\in\mathbb{R}^{m\times n}$
\begin{displaymath}
\begin{split}
\|X+Y\|_{\textup{BiN}}&\leq \alpha\|X+Y\|_{*}\\
&\leq \alpha(\|X\|_{*}+\|Y\|_{*})\\
&\leq \alpha(\|X\|_{\textup{BiN}}+\|Y\|_{\textup{BiN}}).
\end{split}
\end{displaymath}

$\forall X\in\mathbb{R}^{m\times n}$ and $X=UV^{T}$, we have
\begin{displaymath}
\|X\|_{\textup{BiN}}=\min_{X=UV^{T}}\|U\|_{*}\|V\|_{*}\geq0.
\end{displaymath}
In addition, if $\|X\|_{\textup{BiN}}=0$, we have $\|X\|_{*}\leq\|X\|_{\textup{BiN}}=0$, i.e., $\|X\|_{*}=0$.
Hence, $X=0$. In short, the bi-nuclear norm $\|\cdot\|_{\textup{BiN}}$ is a quasi-norm.
\end{proof}

\emph{\textbf{Proof of Theorem 2:}}

\begin{proof}
To prove this theorem, we use the same notations as in the proof of Theorem 1, for instance, $X\!=\!L_{X}\Sigma_{X}R^{T}_{X}$, $U\!=\!L_{U}\Sigma_{U}R^{T}_{U}$ and $V\!=\!L_{V}\Sigma_{V}R^{T}_{V}$ denote the SVDs of $X$, $U$ and $V$, respectively. By Lemma 2, and $p=1/2$, we have
\begin{equation*}
\begin{split}
\|X\|_{S_{1/2}}&\leq\left(\textup{Tr}^{1/2}(O_{2}\Sigma_{X}O^{T}_{2})\right)^{2}=\left(\textup{Tr}^{1/2}(O_{2}O^{T}_{1}\Sigma_{U}O_{4}\Sigma_{V})\right)^{2}=\left(\textup{Tr}^{1/2}(O_{3}\Sigma_{U}O_{4}\Sigma_{V})\right)^{2}\\
&=\left(
\sum^{d}_{i=1}\sqrt{\sum^{d}_{j=1}\tau_{j}(O_{3})_{ij}(O_{4})_{ji}\varrho_{i}}\right )^{2}=\left(
\sum^{d}_{i=1}\sqrt{\varrho_{i}\sum^{d}_{j=1}\tau_{j}(O_{3})_{ij}(O_{4})_{ji}}\right )^{2}\\
&^{a}\!\!\!\leq\sum^{d}_{i=1}\varrho_{i}\sum^{d}_{i=1}\sum^{d}_{j=1}\tau_{j}(O_{3})_{ij}(O_{4})_{ji}\\
&^{b}\!\!\!\leq\sum^{d}_{i=1}\varrho_{i}\sum^{d}_{i=1}\sum^{d}_{j=1}\frac{(O_{3})^{2}_{ij}\tau_{j}+(O_{4})^{2}_{ji}\tau_{j}}{2}\\
&^{c}\!\!\!\leq\sum^{d}_{i=1}\varrho_{i}\sum^{d}_{j=1}\tau_{j}\\
&=\|U\|_{*}\|V\|_{*}\leq\left(\frac{\|U\|_{*}+\|V\|_{*}}{2}\right)^{2},
\end{split}
\end{equation*}
where the inequality $^{a}\!\!\!\leq$ holds due to the Cauchy$-$Schwartz inequality, the inequality $^{b}\!\!\!\leq$ follows from the basic inequality $xy\leq \frac{x^{2}+y^{2}}{2}$ for any real numbers $x$ and $y$, and the inequality $^{c}\!\!\!\leq$ relies on the facts $\sum_{i}(O_{3})^{2}_{ij}=1$ and $\sum_{i}(O_{4})^{2}_{ji}\leq 1$. Thus, we have
\begin{equation*}
\|X\|_{S_{1/2}}\leq\|U\|_{*}\|V\|_{*}\leq\left(\frac{\|U\|_{*}+\|V\|_{*}}{2}\right)^{2}.
\end{equation*}
On the other hand, let $U_{\star}=L_{X}\Sigma^{1/2}_{X}$ and $V_{\star}=R_{X}\Sigma^{1/2}_{X}$, then we have $X=U_{\star}V^{T}_{\star}$ and
\begin{equation*}
\|X\|_{S_{1/2}}=\left(\textup{Tr}^{1/2}(\Sigma_{X})\right)^{2}=\left(\frac{\|L_{X}\Sigma^{1/2}\|_{*}+\|R_{X}\Sigma^{1/2}\|_{*}}{2}\right)^{2}=\left(\frac{\|U_{\star}\|_{*}+\|V_{\star}\|_{*}}{2}\right)^{2}. \end{equation*}
Therefore, under the constraint $X=UV^{T}$, we have
\begin{displaymath}
\begin{split}
\|X\|_{S_{1/2}}=\min_{X=UV^{T}}\|U\|_{*}\|V\|_{*}=\min_{X=UV^{T}}\left(\frac{\|U\|_{*}+\|V\|_{*}}{2}\right)^{2}=\min_{X=UV^{T}}\frac{\|U\|^{2}_{*}+\|V\|^{2}_{*}}{2}=\|X\|_{\textup{BiN}}.
\end{split}
\end{displaymath}
This completes the proof.
\end{proof}
~\\

\section{Proof of Property 3:}
\begin{proof}
The proof involves some following properties of the $\ell_{p}$ quasi-norm, which must be recalled. For any vector $x$ in $\mathbb{R}^{n}$ and $0<p_{2}\leq p_{1}\leq1$, we have
\begin{displaymath}
\|x\|_{\ell_{1}}\leq\|x\|_{\ell_{p_{1}}},\quad \|x\|_{\ell_{p_{1}}}\leq \|x\|_{\ell_{p_{2}}}\leq n^{1/p_{2}-1/p_{1}}\|x\|_{\ell_{p_{1}}}.
\end{displaymath}
Suppose $X\!\in\!\mathbb{R}^{m\times n}$ is of rank $r$, and denote its SVD by $X=U_{m\times r}\Sigma_{r\times r}V^{T}_{n\times r}$. By Theorems 1 and 2, and the properties of the $\ell_{p}$ quasi-norm, we have
\begin{displaymath}
\begin{split}
\|X\|_{\ast}&=\|\textrm{diag}(\Sigma_{r\times r})\|_{\ell_{1}}\leq \|\textrm{diag}(\Sigma_{r\times r})\|_{\ell_{2/3}}=\|X\|_{\textup{F/N}}\leq \sqrt{r}\|X\|_{\ast},\\
\|X\|_{\ast}&=\|\textrm{diag}(\Sigma_{r\times r})\|_{\ell_{1}}\leq \|\textrm{diag}(\Sigma_{r\times r})\|_{\ell_{1/2}}=\|X\|_{\textup{BiN}}\leq r\|X\|_{\ast}.
\end{split}
\end{displaymath}
Similarly,
\begin{displaymath}
\|X\|_{\textup{F/N}}=\|\textrm{diag}(\Sigma_{r\times r})\|_{\ell_{2/3}}\leq\|\textrm{diag}(\Sigma_{r\times r})\|_{\ell_{1/2}}=\|X\|_{\textup{BiN}}.
\end{displaymath}
This completes the proof.
\end{proof}
~\\

\section{Proof of Theorem 3:}
In this paper, the proposed algorithms are based on the proximal alternating linearized minimization (PALM) method for solving the following non-convex problem:
\begin{equation}
\min_{x,y}\Psi(x,y):=F(x)+G(y)+H(x,y),
\end{equation}
where $F(x)$ and $G(y)$ are proper lower semi-continuous functions, and $H(x,y)$ is a smooth function with Lipschitz continuous gradients on any bounded set.

In Algorithm 1, we state that our algorithm alternates between two blocks of variables, $U$ and $V$. We establish the global convergence of Algorithm 1 by transforming the problem (3) into a standard form (11), and show that the transformed problem satisfies the condition needed to establish the convergence. First, the minimization problem (3) can be expressed in the form of (11) by setting
\begin{equation*}
\begin{cases}
F(U):=\frac{2\lambda}{3}\|U\|_{*};\\
G(V):=\frac{\lambda}{3}\|V\|^{2}_{F};\\
H(U,V):=\frac{1}{2}\|\mathcal{P}_{\Omega}(UV^{T})-\mathcal{P}_{\Omega}(D)\|^{2}_{F}.
\end{cases}
\end{equation*}

The conditions for global convergence of the PALM algorithm proposed in~\citep{bolte:palm} are shown in the following lemma.
\begin{lemma}
Let $\{(x_{k},y_{k})\}$ be a sequence generated by the PALM algorithm. This sequence converges to a critical point of (11), if the following conditions hold:
\begin{enumerate}
  \item $\Psi(x,y)$ is a Kurdyka-{\L}ojasiewicz (KL) function;
  \item $\nabla H(x,y)$ has Lipschitz constant on any bounded set;
  \item $\{(x_{k},y_{k})\}$ is a bounded sequence.
\end{enumerate}
\end{lemma}

As stated in the above lemma, the first condition requires that the objective function satisfies the KL property (For more details, see~\citep{bolte:palm}). It is known that any proper closed semi-algebraic function is a KL function as such a function satisfies the KL property for all points in $\textup{dom}f$ with $\varphi(s)=cs^{1-\theta}$ for some $\theta\in[0,1)$ and some $c>0$. Therefore, we first give the following definitions of semi-algebraic sets and functions, and then prove that the proposed problem (3) is also semi-algebraic.

\begin{definition}[\citet{bolte:palm}]
A subset $S\subset\mathbb{R}^{n}$ is a real semi-algebraic set if there exists a finite number of real polynomial functions $\phi_{ij},\,\psi_{ij}:\mathbb{R}^{n}\rightarrow \mathbb{R}$ such that
\begin{equation*}
S=\bigcup_{j}\bigcap_{i} \left\{u\in \mathbb{R}^{n}:\phi_{ij}(u)=0,\;\psi_{ij}(u)<0\right\}.
\end{equation*}
Moreover, a function $h(u)$ is called semi-algebraic if its graph $\{(u,t)\in \mathbb{R}^{n+1}: h(u)=t\}$ is a semi-algebraic set.
\end{definition}

Semi-algebraic sets are stable under the operations of finite union, finite intersections, complementation and Cartesian product. The following are the semi-algebraic functions or the property of semi-algebraic functions used below:
\begin{itemize}
\item Real polynomial functions.
\item Finite sums and product of semi-algebraic functions.
\item Composition of semi-algebraic functions.
\end{itemize}

\begin{lemma}
Each term in the proposed problem (3) is a semi-algebraic function, and thus the function (3) is also semi-algebraic.
\end{lemma}

\begin{proof}
It is easy to notice that the sets $\mathcal{U}=\{U\in \mathbb{R}^{m\times d}:\|U\|_{\infty}\leq D_{1}\}$ and $\mathcal{V}=\{V\in \mathbb{R}^{n\times d}:\|V\|_{\infty}\leq D_{2}\}$ are both semi-algebraic sets, where $D_{1}$ and $D_{2}$ denote two pre-defined upper-bounds for all entries of $U$ and $V$, respectively.

For the first term $F(U)=\frac{2\lambda}{3}\|U\|_{*}$. According to~\citep{bolte:palm}, we can know that the $\ell_{1}$-norm is a semi-algebraic function. Since the nuclear norm is equivalent to the $\ell_{1}$-norm on singular values of the associated matrix, it is natural that the nuclear norm is also semi-algebraic.

For both terms $G(V)=\frac{\lambda}{3}\|V\|^{2}_{F}$ and $H(U,V)=\frac{1}{2}\|\mathcal{P}_{\Omega}(UV^{T}\!-\!D)\|^{2}_{F}$, they are real polynomial functions, and thus are semi-algebraic functions~\citep{bolte:palm}. Therefore, (3) is semi-algebraic due to the fact that a finite sum of semi-algebraic functions is also semi-algebraic.
\end{proof}

For the second condition in Lemma 3, $H(U,V)=\frac{1}{2}\|\mathcal{P}_{\Omega}(UV^{T}\!-\!D)\|^{2}_{F}$ is a smooth polynomial function, and $\nabla H(U,V)=([\mathcal{P}_{\Omega}(UV^{T}\!-\!D)]V,\,[\mathcal{P}_{\Omega}(UV^{T}\!-\!D)]^{T}U)$. It is natural that $\nabla H(U,V)$ has Lipschitz constant on any bounded set.

For the final condition in Lemma 3, $U_{k}\in\mathcal{U}$ and $V_{k}\in\mathcal{V}$ for any $k=1,2,\ldots$, which implies the sequence $\{(U_{k},V_{k})\}$ is bounded.

In short, we can know that three similar conditions as in Lemma 3 hold for Algorithm 1. According to the above discussion, it is clear that the problem (4) is also a semi-algebraic function. In other words, another proposed algorithm shares the same convergence property as in Theorem 3.~\\

\section{Proof of Theorem 5:}
According to Theorem 3, we can know that $(\widehat{U},\widehat{V})$ is a critical point of the problem (3). To prove Theorem 5, we first give the following lemmas.

\begin{lemma}[\citet{lin:almm}]
Let $\mathcal{H}$ be a real Hilbert space endowed with an inner product $\langle\cdot,\cdot\rangle$ and an associated norm $\|\!\cdot\!\|$, and $y\in\partial\|x\|$, where $\partial\|\!\cdot\!\|$ denotes the subgradient of the norm. Then $\|y\|^{*}=1$ if $x\neq 0$, and $\| y\|^{*}\leq 1$ if $x=0$, where $\|\!\cdot\!\|^{*}$ is the dual norm of $\|\!\cdot\!\|$. For instance, the dual norm of the nuclear norm is the spectral norm, $\|\!\cdot\!\|_{2}$, i.e., the largest singular value.
\end{lemma}

\begin{lemma}[\citet{wang:mfcf}]
Let $\mathcal{L}(X)=\frac{1}{\sqrt{mn}}\|X-\widehat{X}\|_{F}$ and $\hat{\mathcal{L}}(X)=\frac{1}{\sqrt{|\Omega|}}\|\mathcal{P}_{\Omega}(X-\widehat{X})\|_{F}$ be the actual and empirical loss function respectively, where $X,\;\widehat{X}\in \mathbb{R}^{m\times n}\;(m\geq n)$. Furthermore, assume entry-wise constraint $\max_{i,j}|X_{ij}|\leq \beta_{1}$. Then for all rank-$r$ matrices $X$, with probability greater than $1-2\exp(-m)$, there exists a fixed constant $C$ such that
\begin{displaymath}
\sup_{X\in S_{r}}|\hat{\mathcal{L}}(X)-\mathcal {L}(X)|\leq C\beta_{1} \left(\frac{mr\log(m)}{|\Omega|}\right)^{1/4},
\end{displaymath}
where $S_{r}=\{X\in \mathbb{R}^{m\times n}: rank(X)\leq r, \|X\|_{F}\leq \sqrt{mn}\beta_{1}\}$.
\end{lemma}

\subsection*{Proof of Theorem 5:}
\begin{proof}
By $C_{2}=\|\mathcal{P}_{\Omega}(D-\widehat{U}\widehat{V}^{T})\widehat{V}\|_{F}/\|\mathcal{P}_{\Omega}(D-\widehat{U}\widehat{V}^{T})\|_{F}$, we have
\begin{equation*}
\begin{split}
&\frac{\|D-\widehat{U}\widehat{V}^{T}\|_{F}}{\sqrt{mn}}\\
\leq&\left|\frac{\|D-\widehat{U}\widehat{V}^{T}\|_{F}}{\sqrt{mn}}-\frac{\|\mathcal{P}_{\Omega}(D-\widehat{U}\widehat{V}^{T})\widehat{V}\|_{F}}{C_{2}\sqrt{|\Omega|}}\right|+\frac{\|\mathcal{P}_{\Omega}(D-\widehat{U}\widehat{V}^{T})\widehat{V}\|_{F}}{C_{2}\sqrt{|\Omega|}}\\
=&\left|\frac{\|D-\widehat{U}\widehat{V}^{T}\|_{F}}{\sqrt{mn}}-\frac{\|\mathcal{P}_{\Omega}(D-\widehat{U}\widehat{V}^{T})\|_{F}}{\sqrt{|\Omega|}}\right|+\frac{\|\mathcal{P}_{\Omega}(D-\widehat{U}\widehat{V}^{T})\widehat{V})\|_{F}}{C_{2}\sqrt{|\Omega|}}.
\end{split}
\end{equation*}
Let $\tau(\Omega):=\left|\frac{1}{\sqrt{mn}}\|D-\widehat{U}\widehat{V}^{T}\|_{F}-\frac{1}{\sqrt{|\Omega|}}\|\mathcal{P}_{\Omega}(D-\widehat{U}\widehat{V}^{T})\|_{F}\right|$, then we need to bound $\tau(\Omega)$. Since rank$(\widehat{U}\widehat{V}^{T})\leq d$ and $\widehat{U}\widehat{V}^{T}\in S_{d}$, and according to Lemma 6, then with probability greater than $1-2\exp(-m)$, then there exists a fixed constant $C_{1}$ such that
\begin{equation}
\begin{split}
\sup_{\hat{U}\hat{V}^{T}\in S_{d}}\tau(\Omega)=&\left|\frac{\|\widehat{U}\widehat{V}^{T}-D\|_{F}}{\sqrt{mn}}-\frac{\|\mathcal{P}_{\Omega}(\widehat{U}\widehat{V}^{T})-\mathcal{P}_{\Omega}(D)\|_{F}}{\sqrt{|\Omega|}}\right|\\
\leq & C_{1} \beta \left(\frac{md\log(m)}{|\Omega|}\right)^{\frac{1}{4}}.
\end{split}
\end{equation}

We also need to bound $\|\mathcal{P}_{\Omega}(\widehat{U}\widehat{V}^{T}-D)\widehat{V}\|_{F}$. Given $\widehat{V}$, the optimization problem with respect to $U$ is formulated as follows:
\begin{equation}
\min_{U} \frac{2\lambda\|U\|_{*}}{3}+\frac{1}{2}\|\mathcal{P}_{\Omega}(U\widehat{V}^{T})-\mathcal{P}_{\Omega}(D)\|^{2}_{F}.
\end{equation}

Since $(\widehat{U}, \widehat{V})$ is a critical point of the problem (3), the first-order optimal condition of the problem (13) is written as follows:
\begin{equation}
\mathcal{P}_{\Omega}(D-\widehat{U}\widehat{V}^{T})\widehat{V}\in \frac{2\lambda}{3}\partial\|\widehat{U}\|_{*}.
\end{equation}
Using Lemma 5, we obtain
\begin{displaymath}
\|\mathcal{P}_{\Omega}(\widehat{U}\widehat{V}^{T}-D)\widehat{V}\|_{2}\leq \frac{2\lambda}{3},
\end{displaymath}
where $\|\!\cdot\!\|_{2}$ is the spectral norm. According to $\textrm{rank}(\mathcal{P}_{\Omega}(\widehat{U}\widehat{V}^{T}-D)\widehat{V})\leq d$, we have
\begin{equation}
\|\mathcal{P}_{\Omega}(\widehat{U}\widehat{V}^{T}-D)\widehat{V}\|_{F}\leq \sqrt{d}\|\mathcal{P}_{\Omega}(\widehat{U}\widehat{V}^{T}-D)\widehat{V}\|_{2} \leq \frac{2\sqrt{d}\lambda}{3}.
\end{equation}
By (12) and (15), we have
\begin{equation*}
\begin{split}
\frac{\|Z-\widehat{U}\widehat{V}^{T}\|_{F}}{\sqrt{mn}}\leq & \frac{\|E\|_{F}}{\sqrt{mn}}+\frac{\|D-\widehat{U}\widehat{V}^{T}\|_{F}}{\sqrt{mn}}\\
\leq & \frac{\|E\|_{F}}{\sqrt{mn}}+\tau(\Omega)+\frac{\|\mathcal{P}_{\Omega}(D-\widehat{U}\widehat{V}^{T})\widehat{V}\|_{F}}{C_{2}\sqrt{|\Omega|}}\\
\leq & \frac{\|E\|_{F}}{\sqrt{mn}}+C_{1}\beta\left(\frac{md\log(m)}{|\Omega|}\right)^{\frac{1}{4}}+\frac{2\sqrt{d}\lambda}{3C_{2}\sqrt{|\Omega|}}.
\end{split}
\end{equation*}
This completes the proof.
\end{proof}

\setcounter{algorithm}{1}

\begin{algorithm}[t]
\caption{Solving (4) via PALM}
\label{alg:Framwork2}
\renewcommand{\algorithmicrequire}{\textbf{Input:}}
\renewcommand{\algorithmicensure}{\textbf{Initialize:}}
\renewcommand{\algorithmicoutput}{\textbf{Output:}}
\begin{algorithmic}[1]
\REQUIRE $\mathcal{P}_{\Omega}(D)$, the given rank $d$ and $\lambda$.
\ENSURE $U_{0}$, $V_{0}$, $\varepsilon$ and $k=0$.\\
\WHILE {not converged}
\STATE {Update $l^{g}_{k+1}$ and $U_{k+1}$ by\\
\begin{displaymath}
l^{g}_{k+1}=\|V_{k}\|^{2}_{2}\;\; \textup{and}\;\;U_{k+1}=\mathop{\arg\min}_{U}\frac{\lambda}{2}\|U\|_{*}+\frac{l^{g}_{k+1}}{2}\|U-U_{k}+\frac{\nabla\! g_{k}(U_{k})}{l^{g}_{k+1}}\|^{2}_{F}.
\end{displaymath} }

\STATE {Update $l^{h}_{k+1}$ and $V_{k+1}$ by\\
\begin{displaymath}
l^{h}_{k+1}=\|U_{k+1}\|^{2}_{2}\;\; \textup{and}\;\;V_{k+1}=\mathop{\arg\min}_{V}\frac{\lambda}{2}\|V\|_{*}+\frac{l^{h}_{k+1}}{2}\|V-V_{k}+\frac{\nabla\! h_{k}(V_{k})}{l^{h}_{k+1}}\|^{2}_{F}.
\end{displaymath} }
\STATE {Check the convergence condition, $\max\{\|U_{k+1}\!-\!U_{k}\|_{F}, \|V_{k+1}\!-\!V_{k}\|_{F}\}<\varepsilon$.}
\ENDWHILE
\OUTPUT $U_{k+1}$, $V_{k+1}$.
\end{algorithmic}
\end{algorithm}

\subsection*{Lower bound on $C_{2}$}
Finally, we also discuss the lower boundedness of $C_{2}$, that is, it is lower bounded by a positive constant. By the characterization of the subdifferentials of norms, we have
\begin{equation}
\partial\|X\|_{*}=\left\{Y\,|\,\langle Y,\,X\rangle=\|X\|_{*},\,\|Y\|_{2}\leq 1\right\}.
\end{equation}
Let $Q=\mathcal{P}_{\Omega}(D-\widehat{U}\widehat{V}^{T})\widehat{V}$, and by (14), we have $Q\in \frac{2\lambda}{3}\partial\|\widehat{U}\|_{*}$. By (16), we have
\begin{displaymath}
\left\langle \frac{3}{2\lambda}Q, \,\widehat{U}\right\rangle=\|\widehat{U}\|_{*}.
\end{displaymath}

Note that $\|X\|_{*}\geq \|X\|_{F}$ and $\langle X,Y\rangle\leq\|X\|_{F}\|Y\|_{F}$ for any matrices $X$ and $Y$ of the same size.
\begin{displaymath}
\frac{3}{2\lambda}\|Q\|_{F}\|\widehat{U}\|_{F}\geq\left\langle \frac{3}{2\lambda}Q, \,\widehat{U}\right\rangle=\|\widehat{U}\|_{*}\geq \|\widehat{U}\|_{F}.
\end{displaymath}
Since $\|\widehat{U}\|_{F}>0$ and $\lambda\neq0$, thus we obtain
\begin{displaymath}
\|\mathcal{P}_{\Omega}(D-\widehat{U}\widehat{V}^{T})\widehat{V}\|_{F}=\|Q\|_{F}\geq \frac{2\lambda}{3}.
\end{displaymath}

$\widehat{U}$ is the optimal solution of the problem (13) with the given matrix $\widehat{V}$, then
\begin{equation*}
\|\mathcal{P}_{\Omega}(D-\widehat{U}\widehat{V}^{T})\|^{2}_{F}<\|\mathcal{P}_{\Omega}(D-\widehat{U}\widehat{V}^{T})\|^{2}_{F}+\frac{2\lambda}{3}\|\widehat{U}\|_{*}\leq\|\mathcal{P}_{\Omega}(D)\|^{2}_{F}=\gamma,
\end{equation*}
where $\gamma>0$ is a constant. Thus,
\begin{displaymath}
C_{2}=\frac{\|\mathcal{P}_{\Omega}(D-\widehat{U}\widehat{V}^{T})\widehat{V}\|_{F}}{\|\mathcal{P}_{\Omega}(D-\widehat{U}\widehat{V}^{T})\|_{F}}>\frac{2\lambda}{3\sqrt{\gamma}}>0.
\end{displaymath}
~\\
%

\section{PALM Algorithm for Solving (4)}
We present an efficient proximal alternating linearized minimization (PALM) algorithm for solving the bi-nuclear quasi-norm regularized matrix completion problem (4), as outlined in \textbf{Algorithm 2}. Moreover, Algorithm 2 shares the same convergence property as in Theorems 3 and 4, and has the similar theoretical recovery guarantee as in Theorem 5.

\bibliography{sigproc}
\bibliographystyle{aaai}